\documentclass[psamsfonts]{amsart}
\usepackage[utf8]{inputenc}
\usepackage{amsmath,amssymb,amsthm,color,enumerate,comment,centernot,enumitem,url,cite, mathrsfs,frcursive,float,stmaryrd,colortbl, lmodern, hyperref}
\usepackage[utf8]{inputenc}
\usepackage{cleveref}
\usepackage{pdfpages} 
\usepackage{pifont}
\usepackage{tikz-cd}
\usepackage{tikz}
\usetikzlibrary{matrix,arrows,decorations.pathmorphing}
\usepackage{graphicx,relsize}
\usepackage{mathtools, enumitem}
\usepackage{array}
\usepackage[margin=1in]{geometry}
\usepackage[utf8]{inputenc}
\usepackage{indentfirst}
\usepackage{amsfonts}
\usepackage{color}
\usepackage{tikz-cd}
\usepackage{setspace}
\usepackage{circuitikz}
\usepackage{mathtools}
\usepackage{mathrsfs}
\usepackage{nicematrix}
\usepackage{lipsum}
\DeclareMathOperator{\ord}{ord}
\DeclareMathOperator{\vol}{Vol}
\DeclareMathOperator{\sgn}{sgn}

\DeclareMathOperator{\tr}{Tr}

\newtheorem{thm}{Theorem}[section]
\newtheorem{cor}[thm]{Corollary}
\newtheorem{prop}[thm]{Proposition}
\newtheorem{lem}[thm]{Lemma}

\newtheorem{Assumption}[thm]{Assumption}

\theoremstyle{definition}
\newtheorem{defn}[thm]{Definition}

\theoremstyle{remark}
\newtheorem{rem}[thm]{Remark}

\numberwithin{equation}{section}

\title{Newton polygons for certain two variable exponential sums}

\author{Bolun Wei}
\address{Institute for math \& AI\\ Wuhan University}
\email[Bolun Wei]{bolunwei@whu.edu.cn}

\begin{document}

\begin{abstract}
Let $f_{t}(x,y)=x^{n}+y+\frac{t}{xy}$ be a Laurent polynomial over $\mathbb{F}_{q}$ with $t$ a parameter. This paper studies the Newton polygon for the $L$-function $L(f_{t},T)$ of toric exponential sums attached to $f_{t}$ over a finite field with characteristic $p$. The explicit Newton polygon is obtained by systematically using Dwork's $\theta_{\infty}$-splitting function with an appropriate choice of basis for cohomology following the method of\cite{AS02}. Our result provides a non-trivial explicit Newton polygon for a non-ordinary family of more than one variable with asymptotical behavior, which gives an evidence of Wan's limit conjecture.
\\\\
\end{abstract}

\maketitle
\tableofcontents

\section{Introduction}
Let $\mathbb{F}_{q}$ be a finite field with $q$ elements of characteristic $p$, let $\zeta_{p}$ be a primitive $p^{\mathrm{th}}$ root of unity in the field of complex numbers. Let $\mathbb{F}_{q^{k}}$ be the finite extension of $\mathbb{F}_{q}$ with degree $k$. For a Laurent polynomial $f\in\mathbb{F}_{q}[x_{1}^{\pm},\cdot\cdot\cdot,x_{m}^{\pm}]$, the toric exponential sum attached to $f$ is defined as 
$$S^{*}_{k}(f)=\sum_{x_{i}\in\mathbb{F}_{q^{k}}^{*}}\zeta_{p}^{\tr_{k}f(x_{1},\cdot\cdot\cdot,x_{m})}$$
where $\mathbb{F}_{q^{k}}^{*}$ denotes the set of non-zero elements in $\mathbb{F}_{q^{k}}$ and $\tr_{k}$ is the trace map from $\mathbb{F}_{q^{k}}$ to $\mathbb{F}_{p}$. By a well-known theorem of Dwork-Bombieri-Grothendieck, the $L$-function is a rational function:
\begin{equation} \label{L}
L(f,T)=\exp(\sum_{k=1}^{\infty}S^{*}_{k}(f)\frac{T^{k}}{k})=\frac{\prod_{i=1}^{d_{1}}(1-\alpha_{i}T)}{\prod_{j=1}^{d_{2}}(1-\beta_{j}T)}    
\end{equation}
where the finitely many numbers $\alpha_{i}$ and $\beta_{j}$ are non-zero algebraic integers. Equivalently we have
$$S^{*}_{k}(f)=\sum_{j=1}^{d_{2}}\beta_{j}^{k}-\sum_{i=1}^{d_{1}}\alpha_{i}^{k}.$$
\par Thus, the study of such $L$-functions is reduced to understanding the reciprocal zeros $\alpha_{i}$ and the reciprocal poles $\beta_{j}$. Without any restriction on $f$, Deligne\cite{Del02} gives some general information about the nature of the roots and poles. For the complex absolute value, we have 
$$|\alpha_{i}|=\sqrt{q^{u_{i}}},\ |\beta_{j}|=\sqrt{q^{v_{i}}},\ u_{i},v_{j}\in\mathbb{Z}\cap[0, 2m],$$
and each $\alpha_{i}$, $\beta_{j}$ and their Galois conjugates over $\mathbb{Q}$ have the same complex absolute value. For a prime $\ell$, denote $\mathbb{Q}_{\ell}$ the field of $\ell$-adic numbers. We tacitly fix an embedding of $\overline{\mathbb{Q}}$ into $\overline{\mathbb{Q}}_{\ell}$, an algebraic closure of $\mathbb{Q}_{\ell}$. When $\ell\neq p$, every $\alpha_{i}$ and $\beta_{j}$ are $\ell$-adic units:
$$|\alpha_{i}|_{\ell}=1,\ |\beta_{j}|_{\ell}=1.$$
When $\ell=p$, we have
$$|\alpha_{i}|_{p}=q^{-r_{i}},\ |\beta_{j}|_{p}=q^{-s_{j}},\ \mathrm{for}\ \mathrm{some}\ r_{i},s_{j}\in\mathbb{Q}\cap[0,m]$$
where we normalize the $p$-adic absolute value such that $|q|_{p}=q^{-1}$. The study of $L$-functions of toric exponential sums is then to study the arithmetic invariants $\{d_{1}, d_{2}, u_{i}, v_{j}, r_{i}, s_{j}\}$. The $u_{i}$ and $v_{j}$ are called the \textbf{weights}
of $\alpha_{i}$ and $\beta_{j}$, and $r_{i}$ and $s_{j}$ are called the \textbf{slopes} of $\alpha_{i}$ and $\beta_{j}$.
\par To get more information about the weights and slopes, we impose a smooth condition on the Laurent polynomial $f$. Write:
$$f(x)=\sum_{w\in\mathbb{Z}^{n}}a_{w}x^{w},\ a_{w}\in\mathbb{F}_{q}$$
where only finitely many $a_{w}$ are non-zero. Here $w=(w_{1},\cdot\cdot\cdot,w_{m})$ is a lattice point in $\mathbb{Z}^{m}$ and $x^{w}$ denotes the monomial $x_{1}^{w_{1}}x_{2}^{w_{2}}\cdot\cdot\cdot x_{m}^{w_{m}}$. We define the \textbf{Newton polytope} of $f$ as
$$\triangle(f)=\mathrm{convex}\ \mathrm{closure}\ \mathrm{of}\ \{0\}\cup\mathrm{Supp}(f)\ \mathrm{in}\ \mathbb{R}^{m}$$
where $\mathrm{Supp}(f)=\{w\in\mathbb{Z}^{m}\mid a_{w}\neq 0\}$. If $\delta$ is a subface of $\triangle(f)$, define the restriction of $f$ to $\delta$ to be the Laurent polynomial 
$$f^{\delta}(x)=\sum_{w\in\delta\cap\mathrm{Supp}(f)}a_{w}x^{w}.$$
\begin{defn} \label{nen-degenerate}
\textit{The Laurent polynomial $f$ is called \textbf{non-degenerate} if for every closed subface $\delta$ of $\triangle(f)$ of arbitrary dimension which does not contain the origin, the Laurent polynomials
$$\frac{\partial f^{\delta}}{\partial x_{1}},\ \frac{\partial f^{\delta}}{\partial x_{2}},\ \cdot\cdot\cdot\ , \ \frac{\partial f^{\delta}}{\partial x_{m}}$$
have no common zero in $(\overline{\mathbb{F}}_{q}^{*})^{m}$.}
\end{defn}
\begin{thm} \label{theorem 1.2}
\normalfont(Adolphson and Sperber\cite{AS01}) \textit{Suppose $f$ is a non-degenerate Laurent polynomial of $m$ variables with coefficients in $\mathbb{F}_{q}$, with $\triangle$ its Newton polytope of dimension $m$, denote $\vol(\triangle)$ the volume of $\triangle$, then we have:\\
(i) $L(f,T)^{(-1)^{m-1}}$ is a polynomial of degree $m!\vol(\triangle).$ \\
(ii) Moreover, if $0$ is an interior point of $\triangle$, then $L(f,T)^{(-1)^{m-1}}$ is pure of weight $m$ (i.e. all reciprocal roots of $L(f,T)^{(-1)^{m-1}}$ have complex absolute value $\sqrt{q^{m}}$).}
\end{thm}
\par This theorem was firstly proved by Adolphson and Sperber\cite{AS01} for almost all primes $p$, later on Denef and Loeser\cite{DL01} proved this for all primes $p$ using the $\ell$-adic method. 
\par Assuming $f$ is non-degenerate, then we may write:
$$ L(f,T)^{(-1)^{m-1}}=\sum_{k=0}^{m!\vol(\triangle)}A_{k}(f)T^{k},\ A_{0}(f)=1,\ A_{k}(f)\in\mathbb{Z}[\zeta_{p}].$$
The \textbf{$q$-adic Newton polygon} $\mathrm{NP}_{q}(f)$ of $L(f,T)^{(-1)^{m-1}}$, is the lower convex hull in $\mathbb{R}^{2}$ of the points 
$$(k,\ord_{q}A_{k}(f)),\ k=0,1,\cdot\cdot\cdot, m!\vol(\triangle),$$
where $\ord_{q}$ the normalized $q$-adic valuation such that $\ord_{q}(q)=1$. It is well-known that the slopes of each line segment in the Newton polygon are the slopes of the reciprocal roots of $L(f,T)^{(-1)^{m-1}}$, and the horizontal length of each line segment is the multiplicity of the reciprocal roots who have the same $q$-adic order. Thus understanding the slopes of the $L$-function turns to the study of the corresponding Newton polygon. 
\par In general determining the exact Newton polygon is a difficult problem  even in low dimensional cases. However, there is a general property that the Newton polygon lies on or above a certain convex hull called the Hodge polygon. We now introduce this combinatorial or topological lower bound. 
\par For $\triangle$ the Newton polytope of $f$, define the cone $\mathrm{Cone}(\triangle)$ to be the union of all rays starting from the origin and passing through $\triangle$, and $M(\triangle)=\mathrm{Cone}(\triangle)\cap\mathbb{Z}^{n}$ the monoid of $\mathbb{Z}$-lattice points lie in the cone. Define the \textbf{weight function} $\omega$ as follow:
\begin{equation} \label{weight}
    \omega: M(\triangle)\longrightarrow\mathbb{R}_{\geq 0}:\ u\mapsto\omega(u):=\mathrm{min}\{c\in\mathbb{R}_{\geq 0}\mid u\in c\triangle\}
\end{equation}
where $c\triangle=\{cx|x\in\triangle\}$ is the dialation of $\triangle$ centered at 0 by a factor $c$.  
\par Note that the image of the weight function is a set of some positive rational numbers. There is a smallest positive integer $D$, called the \textbf{denominator} of $\triangle$, such that the image of $\omega$ lies in $(1/D)\mathbb{Z}_{\geq 0}$. Denote 
$$W_{\triangle}(k)=\#\{u\in M(\triangle)|\omega(u)=\frac{k}{D}\},$$
the number of lattice points in $M(\triangle)$ with weight $k/D$. Define the \textbf{Hodge numbers}
\begin{equation}\label{Hodge number}
    H_{\triangle}(k)=\sum_{i=0}^{k}(-1)^{i}\binom{m}{i}W_{\triangle}(k-iD).
\end{equation}
This number comes from a $p$-adic cohomology space used to compute the $L$-function. $H_{\triangle}(k)$ is a non-negative integer for each $k\in\mathbb{Z}_{\geq 0}$, and for $k>mD$, $H_{\triangle}(k)=0$. Furthermore,
$$\sum_{k=0}^{mD}H_{\triangle}(k)=m!\mathrm{Vol}(\triangle).$$
We define the \textbf{Hodge polygon}, denoted by HP$(\triangle)$ or HP($f$), the lower convex hull in $\mathbb{R}^{2}$ enclosed by the points:
\begin{equation} \label{HP}
    (\sum_{k=0}^{i}H_{\triangle}(k),\frac{1}{D}\sum_{k=0}^{i}kH_{\triangle}(k)),\ i=0,1,2,\cdot\cdot\cdot,mD.
\end{equation}
\par The key result of the Hodge polygon and the Newton polygon is the following theorem:
\begin{thm}
\normalfont(Adolphson and Sperber\cite{AS01}) \textit{For any Laurent polynomial $f$, $\mathrm{NP}_{q}(f)$ lies on or above $\mathrm{HP}(f)$.
The Laurent polynomial $f$ is called \textbf{ordinary} if $\mathrm{NP}_{q}(f)$ equals $\mathrm{HP}(f).$}
\end{thm}
\par Hodge polygons are easier to compute than Newton polygons generally. Thus if a Laurent polynomial is ordinary, we may derive the slopes of reciprocal roots of the $L$-function from the corresponding Hodge polygon. The first example of an ordinary Laurent polynomial family is the Kloosterman sum family
$x+t/x$, studied by Dwork\cite{BD02}. Adolphson and Sperber\cite{AS02}\cite{SS01}\cite{SS02} proved that the hyperkloosterman sum family $x_{1}+\cdot\cdot\cdot+x_{m}+t/(x_{1}\cdot\cdot\cdot x_{m})$ is also an ordinary family whose Newton polygon is the lower convex hull of points $\{(i,i(i-1)/2)\}_{0\leq i\leq n}$. Sperber then studied a generalized hyperkloosterman family $\alpha_{1}x_{1}+\cdot\cdot\cdot+\alpha_{n}x_{n}+tx_{1}^{-a_{1}}x_{2}^{-a_{2}}\cdot\cdot\cdot x_{n}^{-a_{n}}$ in \cite{MR773095} and gave its ordinary condition using Dwork's method. Later on, Bellovin, Garthwaite, Ozman, Pries, Williams, Zhu\cite{MR3204291} obtained the ordinary conditions for $x_{1}^{a_{1}}+\cdot\cdot\cdot+x_{n}^{m_{n}}+x_{1}^{-m_{1}}+\cdot\cdot\cdot x_{n}^{-m_{n}}$ and $x_{1}^{m_{1}}+\cdot\cdot\cdot+x_{n}^{m_{n}}+(x_{1}\cdot\cdot\cdot x_{n})^{-1}$ using Wan's facial decomposition theory\cite{DW01}. More recently, Wang and Yang\cite{WY01} proved that the generalized kloosterman sum family
$f(x_{1},\cdot\cdot\cdot,x_{m})=x_{1}^{a_{1}}+\cdot\cdot\cdot+x_{m}^{a_{m}}+t/(x_{1}^{d_{1}}\cdot\cdot\cdot x_{m}^{d_{m}})$ is ordinary under some congruence condition using the same decomposition theory and Wan's diagonal local theory.
\par However, above examples are either ordinary families, or ordinary under some congruence conditions and explicit Hodge polygons are computed. Newton polygons for non-ordinary families still deserve to be studied. In this paper, we consider the following two variable Laurent polynomial family
$$f_{t}(x,y)=x^{n}+y+\frac{t}{xy},\ \ t\ \mathrm{is}\ \mathrm{a}\ \mathrm{parameter}$$
where $n>1$ is a positive integer. Let
\begin{equation}
    \alpha_{i,j}=i-pj+n\lceil\frac{pj-i}{n}\rceil, \ \ \ i,j\in\mathbb{Z}.   \label{1.5}
\end{equation}
And 
\begin{equation}
    N_{m}=\sum_{i=0}^{m-1}\alpha_{i,\delta(i)}, \ B_{m}=N_{m+1}-N_{m} \label{1.6}
\end{equation}
where $\delta\in S^{0}_{m}=\{\delta\in S_{m}|\sum\limits_{i=0}^{m-1}\alpha_{i,\delta(i)}\ \mathrm{is}\ \mathrm{minimal}\ \mathrm{among}\ \mathrm{all}\ \delta\in S_{m}\}.$
Here elements in $S_{m}$ permutes $0,1,\cdot\cdot\cdot, m-1$ for $0\leq m\leq n+1$. 
\begin{Assumption} \label{assumption 1.4}
\textit{Fix an integer $n>1$, define a Vandermonde-like matrix $$V(x_{0},\cdot\cdot\cdot,x_{m-1})=\left(\begin{array}{ccccc}
1  &  x_{0}^{2} & x_{0}^{2}(x_{0}-1)^{2} & ... & x_{0}^{2}(x_{0}-1)^{2}\cdot\cdot\cdot(x_{0}-m+2)^{2}\\
1  &  x_{1}^{2} & x_{1}^{2}(x_{1}-1)^{2} & ... & x_{1}^{2}(x_{1}-1)^{2}\cdot\cdot\cdot(x_{1}-m+2)^{2}\\
... &  ...& ... &  ...  &  ... \\
1  &  x_{m-1}^{2} & x_{m-1}^{2}(x_{m-1}-1)^{2} & ... & x_{m-1}^{2}(x_{m-1}-1)^{2}\cdot\cdot\cdot(x_{m-1}-m+2)^{2}
\end{array}\right),$$ the determinant of this matrix is non-zero for any set of distinct integers $\{x_{i}\}_{0\leq i\leq m-1}$ with $0\leq x_{i}\leq n-1$. Under this assumption, we denote $$M_{n}(m)=\max\{|\det V(x_{0},\cdot\cdot\cdot,x_{m-1})|:\ x_{0},x_{1},\cdot\cdot\cdot,x_{m-1}\ \mathrm{are}\ \mathrm{distinct},\ \mathrm{and}\ 0\leq x_{i}\leq n-1\}.$$}
\end{Assumption}
Under above assumption, our main results are the following theorems: 
\begin{thm} \label{theorem 1.5}
\textit{For $f_{t}(x,y)=x^{n}+y+\frac{t}{xy}$ where $t\in\mathbb{F}_{p}^{*}$, assumes that \hyperref[assumption 1.4]{assumption 1.4} is satisfied for every $0\leq m\leq n-1$. Then when $p>\max\{M_{n}(0),\cdot\cdot\cdot,M_{n}(n-1),2n^{3}-n^{2}-n+1\}$, the $p$-Newton polygon for $L(f_{t},T)^{-1}$ is the end-to-end join of $2n+1$ line segments of horizontal length 1 with slopes:
$$\{\frac{i}{n}+\frac{(2n+1)B_{i}}{n(p-1)}\}_{0\leq i\leq n}\cup\{\frac{i}{n}-\frac{(2n+1)B_{2n-i}}{n(p-1)}\}_{n+1\leq i\leq 2n.}$$} 
\end{thm}
\par Furtherly if we impose a condition on the base prime $p$, we will obtain the $q$-adic Newton polygon when the parameter $t\in\mathbb{F}_{q}^{*}$ for some $q=p^{a}$:
\begin{Assumption} \label{assumption 1.6}
\textit{For prime $p>n$, it satisfies
$$\ord_{p}[(k-1)!(p-k)!-(-1)^{k}]=1\ \mathrm{for}\ \mathrm{any}\  1\leq k\leq n-1.$$} 
\end{Assumption}
And here is the main theorem for the parameter $t\in\mathbb{F}_{q}^{*}$:
\begin{thm} \label{theorem 1.7}
 \textit{For the family $f_{t}$ with $t\in\mathbb{F}_{q}^{*}$, when the base prime $p>4n^{4}+4n^{3}+3n^{2}+n+1$ satisfies \hyperref[assumption 1.6]{assumption 1.6}, and \hyperref[assumption 1.4]{assumption 1.4} is satisfied for all integer $2\leq m\leq n-1$, then the $q$-adic Newton polygon for $L(f_{t},T)^{-1}$ coincides with the $p$-adic Newton polygon described in \hyperref[theorem 1.5]{theorem 1.5}.}   
\end{thm}
\par As an application of the main results, we compute exact Newton polygons for $n=3,4$: 
\begin{cor} \label{corollary 3.12} Suppose $f_{t}(x,y)=x^{3}+y+\frac{t}{xy}$, then we have: \\
(a) When $p\equiv 1\mod 3$, $f_{t}$ is ordinary. The slope sequence of the Newton polygon is
$$\{0,\frac{1}{3},\frac{2}{3},1,\frac{4}{3},\frac{5}{3},2\},$$
where each ling segment has horizontal length 1. \\
(b) When $p\equiv 2\mod 3$, $t\in\mathbb{F}_{p}^{*}$ and $p>43$, the slope sequence of $\mathrm{NP}_{p}(f_{t})$ is 
$$\{0,\  \frac{1}{3}+\frac{14}{3(p-1)},\ \frac{2}{3}-\frac{14}{3(p-1)},\ 1,\ \frac{4}{3}+\frac{14}{3(p-1)},\ \frac{5}{3}-\frac{14}{3(p-1)},\ 2\}$$
where each line segment has horizontal length 1. \\
(c) When $p\equiv 2\mod 3$, $t\in\mathbb{F}_{q}^{*}$ for some $q=p^{a}$ with $a>1$, $p>463$ and $$\ord_{p}[(p-1)!+1]=\ord_{p}[(p-2)!-1]=1,$$
$\mathrm{NP}_{q}(f_{t})$ coincides with that of case (b).
\end{cor}
\begin{cor} \label{corollary 3.13} Suppose $f_{t}(x,y)=x^{4}+y+\frac{t}{xy}$, then we have: \\
(a) When $p\equiv 1\mod 4$, $f_{t}$ is ordinary. The slope sequence of the Newton polygon is
$$\{0,\frac{1}{4},\frac{1}{2},\frac{3}{4},1,\frac{5}{4},\frac{3}{2},\frac{7}{4},2\},$$
where each ling segment has horizontal length 1. \\
(b) When $p\equiv 3\mod 4$, $t\in\mathbb{F}_{p}^{*}$ and $p>109$, the slope sequence of $\mathrm{NP}_{p}(f_{t})$ is 
$$\{0,\ \frac{1}{4}+\frac{18}{4(p-1)},\ \frac{1}{2},\ \frac{3}{4}-\frac{18}{4(p-1)},\ 1,\ \frac{5}{4}+\frac{18}{4(p-1)},\ \frac{3}{2},\ \frac{7}{4}-\frac{18}{4(p-1)},\ 2\}$$
where each line segment has horizontal length 1. \\
(c) When $p\equiv 3\mod 4$, $t\in\mathbb{F}_{q}^{*}$ for some $q=p^{a}$ with $a>1$, $p>1333$ and $$\ord_{p}[(p-1)!+1]=\ord_{p}[(p-2)!-1]=\ord_{p}[2(p-3)!+1]=1,$$
$\mathrm{NP}_{q}(f_{t})$ coincides with that of case (b).
\end{cor}
\begin{rem}
Numerical calculation shows that \hyperref[assumption 1.4]{assumption 1.4} is true for $n$ large to $10^{6}$, we hope some further combinatoric and linear algebra study can help remove this assumption.   
\end{rem}
\begin{rem}
Notice that the Newton polygon is independent of the choice of the ground field $\mathbb{F}_{q}$, and the Hodge polygon only depends on the combinatorial shape of the Newton polytope. We naturally want to know the behavior of the Newton polygon when the prime $p$ varies and the Laurent polynomial varies in some parameter family. Consider a Laurent polynomial $f$ with coefficients in $\overline{\mathbb{Q}}$, and with Newton polytope $\triangle$. For a prime $p$, fix an embedding $\overline{\mathbb{Q}}\rightarrow\overline{\mathbb{Q}}_{p}$ and view $f$ as a Laurent polynomial with coefficients in $\overline{\mathbb{Q}}_{p}$. Denote $f\ \mathrm{mod}\ p$, the reduction of $f$ with coefficients in the residue field $\mathbb{F}_{q}$ for some $q=p^{a}$. Wan\cite{DW01} conjectured that under some conditions, when the base prime $p$ grows to infinity, the Newton polygon will asymptotically approach to the Hodge polygon:
$$\mathop{\mathrm{lim}}\limits_{p\rightarrow\infty}\mathrm{NP}(f\ \mathrm{mod}\ p)=\mathrm{HP}(\triangle).$$
Zhu\cite{Zhu01}\cite{Zhu02} proved this conjecture for one variable polynomial families. But so far the conjecture still remains widely open.  
\end{rem}
\par Back into our example, we readily compute the Hodge polygon of the family $f_{t}(x,y)=x^{n}+y+\frac{t}{xy}$ to be the end-to-end join $2n+1$ line segments of horizontal length 1 with slopes $\{i/n\}_{0\leq i\leq 2n}$. We therefore give a confirmed answer to Wan's limit conjecture for our family in the following sense:
\begin{cor}
For a fixed $n$ with $t\in\mathbb{Z}\setminus\{0\}$, suppose \hyperref[assumption 1.4]{assumption 1.4} is satisfied, then we have 
$$\mathop{\mathrm{lim}}\limits_{p\rightarrow\infty}\mathrm{NP}_{p}(f_{t}\ \mathrm{mod}\ p)=\mathrm{HP}(\triangle).$$
\end{cor}

\textbf{Acknowledgments.} This paper partly comes from the author's Ph.D. thesis. The author thanks Douglas Haessig for many guidance and encouragement. Also much thanks to Daqing Wan and Steven Sperber for many enlightening conversations through the project.

\section{Dwork cohomology}
\par Through all the paper, $n$ is a fixed positive integer, $\mathbb{F}_{q}$ is a finite field with characteristic $p>2$, $q=p^{a}$, and $p\nmid n$. $\mathbb{Q}_{q}$ the unramified extension of $\mathbb{Q}_{p}$ of degree $a$ and let $\mathbb{Z}_{q}$ be its ring of integers. Fix $\zeta_{p}$ a primitive $p^{\mathrm{th}}$ root of unity in $\overline{\mathbb{Q}}_{p}$. Let $\Omega_{1}=\mathbb{Q}_{p}(\zeta_{p})$, the totally ramified extension for $\mathbb{Q}_{p}$ of degree $p-1$, with ring of integers $\mathcal{O}_{1}=\mathbb{Z}_{p}[\zeta_{p}]$. Denote $\Omega_{0}=\mathbb{Q}_{q}(\zeta_{p})$, with ring of integers $\mathcal{O}_{0}=\mathbb{Z}_{q}[\zeta_{p}]$. Let $\mathbb{C}_{p}$ be the completion of $\Bar{\mathbb{Q}}_{p}$ w.r.t. the $p$-adic norm $|\ \ |_{p}$, then $\mathbb{C}_{p}$ is complete and algebraically closed. For all $t\in\mathbb{F}_{q}^{*}$, the Newton polytope $\triangle$ for our family $f_{t}(x,y)$ is an triangle with 3 vertices $(-1,-1)$, $(n,0)$, $(0,1)$. So Cone($\triangle$) will be all the $\mathbb{R}^{2}$ plane and the monoid $M(\triangle)$ will be all the $\mathbb{Z}$-lattice points in the plane. The corresponding weight function of $f_{t}$ will be
\begin{equation}
\omega:\mathbb{Z}^{2}\longrightarrow\frac{1}{n}\mathbb{Z}_{\geq 0}\ \ (a,b)\mapsto\frac{a}{n}+b+\frac{2n+1}{n}m(a,b) \label{2.1}
\end{equation}
where 
\begin{equation}
m(a,b)=\max\{0,-a,-b\} .  \label{2.2}  
\end{equation}
The weight function satisfies the following property:
\begin{prop}
    Let $\omega$ be the weight function on $M(\triangle)$, then we have:\\
(a) $\omega(u)=0$ if and only if $u=\Vec{0}$ in $\mathbb{R}^{m}$.\\
(b) $\omega(cu)=c\omega(u)$ for any $c\in\mathbb{Z}_{\geq 0}$\\
(c) $\omega(u+v)\leq\omega(u)+\omega(v)$, the equality holds if and only if $u$ and $v$ are co-facial.
\end{prop}
We see that the denominator of $\triangle$ in our family is $n$, volume of the polytope $\mathrm{Vol}(\triangle)=(2n+1)/2$. And the Laurent polynomials in this family are all non-degenerated if $p\nmid n$.
\par Let $\gamma$ be a zero of the power series $\sum_{k=0}^{\infty}\frac{x^{p^{k}}}{p^{k}}$ in $\Omega_{1}$ with $\ord_{p}(\gamma)=\frac{1}{p-1}$. $\mathrm{E}(x)=\exp(\sum_{k=0}^{\infty}\frac{x^{p^{k}}}{p^{k}})$ denotes the Artin-Hasse series, and $\Theta_{\infty}(x)=\mathrm{E}(\gamma x)$ denotes the splitting function for $\gamma$ in Dwork's terminology, this function holds the following properties:
\begin{prop} \label{proposition 2.2}
\normalfont(Dwork\cite{BD03}, §4)   \textit{The splitting function $\Theta_{\infty}(x)=\sum_{k=0}^{\infty}a_{k}x^{k}$ satisfies:\\
(a)  $a_{k}$'s lie in a finite extension of $\mathbb{Q}_{p}$ and $\ord_{p}(a_{k})\geq\frac{k}{p-1}$ for all non-negative integer $k$. In particular, $a_{k}=\frac{\gamma^{k}}{k!}$ and $\ord_{p}(a_{k})=\frac{k}{p-1}$ for $0\leq k\leq p-1$.\\
(b) $\Theta_{\infty}(x)$ converges in the disk $\{x\in\mathbb{C}_{p}|\ord_{p}(x)>-\frac{1}{p-1}\}$.\\
(c) $\Theta_{\infty}(1)$ is a primitive $p^{\mathrm{th}}$ root of unity.\\
(d) If $\Bar{\alpha}\in\mathbb{Q}_{p}$ is a Teichmüller lift of $\alpha\in\mathbb{F}_{q}$ where $q=p^{a}$ for some positive integer $a$ (i.e. $\Bar{\alpha}^{p^{a}}=\Bar{\alpha}$), then we have $$\Theta_{\infty}(1)^{\sum\limits_{k=0}^{a-1}\Bar{\alpha}^{p^{k}}}=\prod\limits_{k=0}^{a-1}\Theta_{\infty}(\Bar{\alpha}^{p^{k}}).$$}
\end{prop}
\par Now let $\Bar{t}\in\mathbb{Q}_{q}$ be a Teichmüller of $t\in\mathbb{F}_{q}^{*}$, and let
\begin{equation}
 F_{t}(x,y)=\Theta_{\infty}(x^{n})\Theta_{\infty}(y)\Theta_{\infty}(\frac{\Bar{t}}{xy})=\sum_{(a,b)\in\mathbb{Z}^{2}}B(a,b)x^{a}y^{b}   \label{2.3} \end{equation}
We have the $p$-adic estimates of the coefficients $B(a,b)$:
\begin{lem} \label{lemma 2.3}
For all $(a,b)\in\mathbb{Z}^{2}$, $\ord_{p}B(a,b)\geq\frac{\omega(a,b)}{p-1}$.
\end{lem}
\textit{Proof.} Expand the coefficients in $F_{t}(x,y)$, we get
$$B(a,b)=\sum_{(k,l,m)\in I(a,b)}\Bar{t}^{k}a_{k}a_{l}a_{m}$$
where $I(a,b)=\{(k,l,m)\in\mathbb{Z}^{3}_{\geq 0}|nl-k=a,\ m-k=b\}$.
\par $\Bar{t}$ is a Teichmüller, then $\Bar{t}\in\mathbb{Z}_{q}$ and $\ord_{p}(\Bar{t})=1$. Apply \hyperref[proposition 2.2]{proposition 2.2} part (a) we get 
$$\ord_{p}B(a,b)\geq\mathop{\mathrm{inf}}\limits_{(k,l,m)\in I(a,b)}\frac{k+l+m}{p-1}.$$
\par For $(k,l,m)\in I(a,b)$, $k=nl-a\geq -a$, $k=m-b\geq -b$, so $k\geq m(a,b)$ where $m(a,b)$ defined in \hyperref[2.2]{(2.2)}. We substitute $l=\frac{a+k}{n}$, $m=b+k$ and use the weight function formula in \hyperref[2.1]{(2.1)} to obtain the estimation
$$\frac{k+l+m}{p-1}=\frac{k+\frac{a+k}{n}+b+k}{p-1}=\frac{\frac{a}{n}+b+\frac{2n+1}{n}k}{p-1}\geq\frac{\frac{a}{n}+b+\frac{2n+1}{n}m(a,b)}{p-1}=\frac{\omega(a,b)}{p-1}.$$
\par In particular, notice that if $k>m(a,b)$, we have $\ord_{p}B(a,b)>\frac{\omega(a,b)}{p-1}.$ $\hfill\square$ \\
\par We now fix $\widetilde{\gamma}$ a root of $x^{n}-\gamma=0$ in $\mathbb{C}_{p}$, note that the ring of integers for $\Omega_{0}(\widetilde{\gamma})$ (resp. $\Omega_{1}(\widetilde{\gamma})$) is $\mathbb{Z}_{q}[\widetilde{\gamma}]$ (resp. $\mathbb{Z}_{p}[\widetilde{\gamma}]$). Then we define a space of $p$-adic functions
\begin{equation}
\mathcal{C}_{0}=\{\sum_{(a,b)\in\mathbb{Z}^{2}}\xi(a,b)\widetilde{\gamma}^{n\omega(a,b)}x^{a}y^{b}|\xi(a,b)\in\mathbb{Z}_{q}[\widetilde{\gamma}],\ |\xi(a,b)|_{p}\rightarrow 0\ \mathrm{as}\ \omega(a,b)\rightarrow \infty\}  
\end{equation}
endowed with the norm
$$|\xi|=\mathop{\mathrm{sup}}\limits_{(a,b)\in\mathbb{Z}^{2}}\{|\xi(a,b)|_{p}\}$$
for $\xi=\sum\limits_{(a,b)\in\mathbb{Z}^{2}}\xi(a,b)\widetilde{\gamma}^{n\omega(a,b)}x^{a}y^{b}\in\mathcal{C}_{0}$. Then $\mathcal{C}_{0}$ is a Banach $\mathbb{Z}_{q}[\widetilde{\gamma}]$-algebra w.r.t. the superior norm.
\par Let $\sigma$ be the Frobenius generator of $Gal(\mathbb{Q}_{q}/\mathbb{Q}_{p})$, then we extend it to $Gal(\Omega_{0} (\widetilde{\gamma})/\Omega_{1}(\widetilde{\gamma}))$ by fixing $\sigma(\widetilde{\gamma})=\widetilde{\gamma}$ and $\sigma(\zeta_{p})=\zeta_{p}$. $\psi_{p}$ be the inverse Frobenius operator acting on $\mathcal{C}_{0}$ by
\begin{equation}
\psi_{p}:\mathcal{C}_{0}\rightarrow\mathcal{C}_{0}\  \sum\limits_{(a,b)\in\mathbb{Z}^{2}}\xi(a,b)\widetilde{\gamma}^{n\omega(a,b)}x^{a}y^{b}\mapsto\sum\limits_{(a,b)\in\mathbb{Z}^{2}}\xi(pa,pb)\widetilde{\gamma}^{n\omega(pa,pb)}x^{a}y^{b}.   
\end{equation}
Define a semi-linear (over $\Omega_{0}(\widetilde{\gamma})$) operator $\alpha_{1}$ by
\begin{equation}
\alpha_{1}=\sigma^{-1}\circ\psi_{p}\circ F_{t}(x,y)    \label{2.6}
\end{equation}
where the composition for $F_{t}(x,y)$ is the multiplication by $F_{t}(x,y)$, $\sigma^{-1}$ acts on the coefficients of the elements in $\mathcal{C}_{0}$. Let $\alpha_{0}=\alpha_{1}^{a}$. Then $\alpha_{0}$ is a completely continuous operator, linear over $\Omega_{0}(\widetilde{\gamma})$ in the sense of\cite{JS01}. So $\alpha_{0}$ has a $p$-adically entire Fredholm determinant, det$(I-T\alpha_{0})$. Let $\delta$ acts on power series via
$$P(T)^{\delta}=\frac{P(T)}{P(qT)}.$$
Together with the Dwork trace formula (\cite{BD01}, lemma 2) 
\begin{equation}
S_{k}(f_{t})=(q^{k}-1)^{2}\mathrm{Tr}(\alpha_{0}^{k})
\end{equation}
and the matrix expression $\mathrm{det}(I-T\alpha_{0})=\mathrm{exp}(-\sum_{k=1}^{\infty}\mathrm{Tr}(\alpha_{0}^{k})\frac{t^{k}}{k})$, we are able to derive the expression of the $L$-function for our family $f_{t}$:
\begin{equation}
L(f_{t},T)^{-1}=\mathrm{det}(I-T\alpha_{0})^{\delta^{2}}. \label{2.8}
\end{equation} 
\par We introduce the cohomology theory to get a cohomological expression of $L$-functions and then compute the Newton polygon. Let $\gamma_{0}=\gamma$, the root of $\sum\limits_{k=0}^{\infty}\frac{x^{p^{k}}}{p^{k}}=0$, for $i\geq 1$, let
$$r_{i}=\sum_{k=0}^{i}\frac{\gamma^{p^{k}}}{p^{k}}=-\sum_{k=i+1}^{\infty}\frac{\gamma^{p^{k}}}{p^{k}}.$$
Use the second description we have 
\begin{equation}
    \ord_{p}(\frac{\gamma_{i}}{\gamma_{0}})=\frac{p^{i+1}-1}{p-1}-(i+1) \label{2.9}
\end{equation}
for all $i\geq 0$. Then we see that $F_{t}(x,y)$ defined in \hyperref[2.3]{(2.3)} can be expressed as 
$$F_{t}(x,y)=\frac{\mathrm{exp}(H_{t}(x,y))}{\mathrm{exp}(H_{t}(x^{p},y^{p}))}$$
where
$$H_{t}(x,y)=\sum_{i=0}^{\infty}\gamma_{i}(x^{np^{i}}+y^{p^{i}}+\frac{\sigma^{i}(\Bar{t})}{x^{p^{i}}y^{p^{i}}}).$$
Here $\Bar{t}$ is a Teichmüller of $t$, so $\sigma(\Bar{t})=\Bar{t}^{p}$. Then we find that the operators $\alpha_{0}$ and $\alpha_{1}$ can be written as
$$\alpha_{1}=\frac{1}{\mathrm{exp}(H_{t}(x,y))}\circ\sigma^{-1}\circ\psi_{p}\circ\mathrm{exp}(H_{t}(x,y)),$$
$$\alpha_{0}=\frac{1}{\mathrm{exp}(H_{t}(x,y))}\circ\psi_{p}^{a}\circ\mathrm{exp}(H_{t}(x,y)).$$
Motivated by this, we define the differential operators on $\mathcal{C}_{0}$ as 
$$D_{x}=\frac{1}{\exp(H_{t}(x,y))}\circ x\frac{\partial}{\partial x}\circ\exp(H_{t}(x,y)),$$
$$D_{y}=\frac{1}{\mathrm{exp}(H_{t}(x,y))}\circ y\frac{\partial}{\partial y}\circ\exp(H_{t}(x,y)).$$
And they can be expressed as
\begin{equation}
D_{x}=x\frac{\partial}{\partial x}+x\frac{\partial H_{t}}{\partial x}=x\frac{\partial}{\partial x}+\sum_{i=0}^{\infty}r_{i}p^{i}(nx^{np^{i}}-\frac{\sigma^{i}(\Bar{t})}{x^{p^{i}}y^{p^{i}}}),
\end{equation}
\begin{equation}
D_{y}=y\frac{\partial}{\partial y}+y\frac{\partial H_{t}}{\partial y}=y\frac{\partial}{\partial y}+\sum_{i=0}^{\infty}r_{i}p^{i}(y^{p^{i}}-\frac{\sigma^{i}(\Bar{t})}{x^{p^{i}}y^{p^{i}}}).    
\end{equation}
We construct the complex $(\Omega_{\mathcal{C}_{0}}^{\bullet},\triangledown(D))$ as in\cite{AS01}
$$\Omega_{\mathcal{C}_{0}}^{0}=\mathcal{C}_{0},\ \Omega_{\mathcal{C}_{0}}^{1}=\mathcal{C}_{0}\frac{dx}{x}\oplus \mathcal{C}_{0}\frac{dy}{y},\ \Omega_{\mathcal{C}_{0}}^{2}=\mathcal{C}_{0}\frac{dx}{x}\wedge\frac{dy}{y},$$
with the boundary map
$$D^{(0)}:\Omega_{\mathcal{C}_{0}}^{0}\rightarrow\Omega_{\mathcal{C}_{0}}^{1},\ \ \xi\mapsto D_{x}(\xi)\frac{dx}{x}+D_{y}(\xi)\frac{dy}{y},$$ $$D^{(1)}:\Omega_{\mathcal{C}_{0}}^{1}\rightarrow\Omega_{\mathcal{C}_{0}}^{2},\ \ \xi_{1}\frac{dx}{x}+\xi_{2}\frac{dy}{y}\mapsto (D_{x}(\xi_{2})-D_{y}(\xi_{1}))\frac{dx}{x}\wedge\frac{dy}{y}.$$
Furthermore, by \cite{BD03} (equation 4.35) we have 
$$\alpha_{1}\circ D_{x}=pD_{x}\circ\alpha_{1},\ \mathrm{and}\ \alpha_{1}\circ D_{y}=pD_{y}\circ\alpha_{1}$$
and therefore
$$\alpha_{0}\circ D_{x}=qD_{x}\circ\alpha_{0},\ \mathrm{and}\ \alpha_{0}\circ D_{y}=qD_{y}\circ\alpha_{0}.$$
Then we can define the Frobenius chain maps
\begin{align*}
&\mathrm{Frob}_{0}^{(0)}:\Omega_{\mathcal{C}_{0}}^{0}\rightarrow\Omega_{\mathcal{C}_{0}}^{0}\ \ \ \xi\mapsto q^{2}\alpha_{0}(\xi),\\
&\mathrm{Frob}_{0}^{(1)}:\Omega_{\mathcal{C}_{0}}^{1}\rightarrow\Omega_{\mathcal{C}_{0}}^{1}\ \ \ \xi_{1}\frac{dx}{x}+\xi_{2}\frac{dy}{y}\mapsto q\alpha_{0}(\xi_{1})\frac{dx}{x}+q\alpha_{0}(\xi_{2})\frac{dy}{y},\ \mathrm{and}\\
&\mathrm{Frob}_{0}^{(2)}:\Omega_{\mathcal{C}_{0}}^{2}\rightarrow\Omega_{\mathcal{C}_{0}}^{2}\ \ \ \xi\frac{dx}{x}\wedge\frac{dy}{y}\mapsto\alpha_{0}(\xi)\frac{dx}{x}\wedge\frac{dy}{y}.    
\end{align*}
With an abuse using of notation, we still denote the maps on the cohomology level as $\mathrm{Frob}_{0}^{\bullet}$, notice all the chain maps are completely continuous operators, $\mathrm{Frob}_{0}^{\bullet}$ are nuclear, therefore we can refine the $L$-function expression in \hyperref[2.8]{(2.8)} as
$$L(f_{t},T)^{-1}=\prod_{i=0}^{2}\mathrm{det}(I-T\mathrm{Frob}_{0}^{(i)}|_{H^{i}(\Omega_{\mathcal{C}_{0}}^{\bullet})})^{(-1)^{i}}$$
where each factor on the right is $p$-adically entire. By Adolphson and Sperber\cite{AS01}, this cohomology is acyclic except $H^{2}(\Omega_{\mathcal{C}_{0}})$ a free $\mathbb{Z}_{q}[\widetilde{\gamma}]$-module of rank $2n+1$ due to the non-degeneracy of $f_{t}(x,y)$. Therefore, the $L$-function for $f_{t}$ can be written as
\begin{equation}
L(f_{t},T)^{-1}=\mathrm{det}(I-T\mathrm{Frob}_{0}^{(2)}|_{H^{2}(\Omega_{\mathcal{C}_{0}}^{\bullet})}),    
\end{equation}
which is a polynomial of degree $2n+1$. The top cohomology $H^{2}(\Omega_{\mathcal{C}_{0}}^{\bullet})\simeq\mathcal{C}_{0}/(D_{x}\mathcal{C}_{0}+ D_{y}\mathcal{C}_{0})$, and $\mathrm{Frob}_{0}^{(2)}$ acts on it as $\alpha_{0}$. We naturally want to find a basis for the top cohomological space and express the explicit matrix w.r.t the basis. To do this, we introduce the reduction cohomology. We define an increasing filtration of $\mathbb{F}_{q}[x^{\pm},y^{\pm}]$ indexed by $i\in\mathbb{Z}_{\geq 0}$ as 
$$\mathrm{Fil}^{i}\mathbb{F}_{q}[x^{\pm},y^{\pm}]=\{\Bar{\xi}=\sum_{(a,b)\in\mathbb{Z}^{2}}\Bar{\xi}(a,b)x^{a}y^{b}|\omega(a,b)\leq\frac{i}{n}\ \mathrm{for}\ \mathrm{all}\ (a,b)\in\mathrm{Supp}(\Bar{\xi})\},$$
if $i<0$ we set $\mathrm{Fil}^{i}\mathbb{F}_{q}[x^{\pm},y^{\pm}]=0$. Let
$$\Bar{S}^{i}=\mathrm{Fil}^{i}\mathbb{F}_{q}[x^{\pm},y^{\pm}]/\mathrm{Fil}^{i-1}\mathbb{F}_{q}[x^{\pm},y^{\pm}].$$
We see $\Bar{S}^{i}\simeq\{\Bar{\xi}=\sum_{(a,b)\in\mathbb{Z}^{2}}\Bar{\xi}(a,b)x^{a}y^{b}|\omega(a,b)=\frac{i}{n}\ \mathrm{for}\ \mathrm{all}\ (a,b)\in\mathrm{Supp}(\Bar{\xi})\}$, and for $i<0$ we set $\Bar{S}^{i}=0$. Let $\Bar{S}$ be the associated graded ring $\mathrm{gr}\mathbb{F}_{q}[x^{\pm},y^{\pm}]=\oplus\Bar{S}^{i}$ where the multiplication is defined as
$$x^{a}y^{b}\cdot x^{c}y^{d}=\left\{\begin{array}{cc}
    x^{a+c}y^{b+d} & \ \ \ \ \  \mathrm{when}\ (a,b),(c,d)\ \mathrm{are}\ \mathrm{cofacial}\ \mathrm{in}\ \triangle, \\
    0 & \ \mathrm{otherwise.}    
\end{array}\right.$$
Then we define a map
\begin{equation}
\mathrm{Pr}:\mathcal{C}_{0}\rightarrow\Bar{S}\ \sum_{(a,b)\in\mathbb{Z}^{2}}\xi(a,b)\widetilde{\gamma}^{n\omega(a,b)}x^{a}y^{b}\mapsto\sum_{(a,b)\in\mathbb{Z}^{2}}\Bar{\xi}(a,b)x^{a}y^{b}    \label{2.13}
\end{equation}
where $\Bar{\xi}(a,b)$ is the reduction of $\xi(a,b)$ in the residue field $\mathbb{F}_{q}$. Pr is a ring homomorphism (\cite{AS01} Lemma 2.10) with $\mathcal{C}_{0}/\widetilde{\gamma}\mathcal{C}_{0}\simeq\Bar{S}$, mapping as a reduction modulo $\widetilde{\gamma}$.
\par By \hyperref[2.9]{(2.9)} $\ord_{p}(r_{i}p^{i})>\frac{p^{i}}{p-1}$ for $i>0$, the higher order terms in $x\frac{\partial H_{t}}{\partial x}$ and $y\frac{\partial H_{t}}{\partial y}$ vanish via the reduction map Pr and only the terms for $i=0$ remains. We have 
$$H_{x}=\mathrm{Pr}(\gamma(nx^{n}-\frac{\Bar{t}}{xy}))=nx^{n}-\frac{t}{xy},\ \mathrm{and}\ H_{y}=\mathrm{Pr}(\gamma(y-\frac{\Bar{t}}{xy}))=y-\frac{t}{xy}.$$
Therefore the reduction differential operator for $D_{x}$, $D_{y}$ mod $\widetilde{\gamma}$ will be
\begin{equation}
\Bar{D}_{x}=x\frac{\partial}{\partial x}+H_{x}, \label{2.14}
\end{equation}
\begin{equation}
\Bar{D}_{y}=y\frac{\partial}{\partial y} +H_{y} .  \label{2.15}
\end{equation}
\par We then construct two complexes on $\Bar{S}$, $(\Omega_{\Bar{S}}^{\bullet},\triangledown(H))$ and $(\Omega_{\Bar{S}}^{\bullet},\triangledown(\Bar{D}))$ as follows. The spaces in both cases are the same:
$$\Omega_{\Bar{S}}^{0}=\Bar{S},\ \Omega_{\Bar{S}}^{1}=\Bar{S}\frac{dx}{x}\oplus \Bar{S}\frac{dy}{y},\ \Omega_{\Bar{S}}^{2}=\Bar{S}\frac{dx}{x}\wedge\frac{dy}{y},$$
where the boundary map for $\triangledown(H)$:
\begin{align*}
&H^{(0)}:\Omega_{\Bar{S}}^{0}\rightarrow\Omega_{\Bar{S}}^{1}\ \ \ \Bar{\xi}\mapsto H_{x}\Bar{\xi}\frac{dx}{x}+H_{y}\Bar{\xi}\frac{dy}{y},\\ &H^{(1)}:\Omega_{\Bar{S}}^{1}\rightarrow\Omega_{\Bar{S}}^{2}\ \ \ \Bar{\xi}_{1}\frac{dx}{x}+\Bar{\xi}_{2}\frac{dy}{y}\mapsto (H_{x}\Bar{\xi}_{2}-H_{y}\Bar{\xi}_{1})\frac{dx}{x}\wedge\frac{dy}{y},   
\end{align*}
and for $\triangledown(\Bar{D})$:
\begin{align*}
&\Bar{D}^{(0)}:\Omega_{\Bar{S}}^{0}\rightarrow\Omega_{\Bar{S}}^{1}\ \ \ \Bar{\xi}\mapsto \Bar{D}_{x}(\Bar{\xi})\frac{dx}{x}+\Bar{D}_{y}(\Bar{\xi})\frac{dy}{y},\\
&\Bar{D}^{(1)}:\Omega_{\Bar{S}}^{1}\rightarrow\Omega_{\Bar{S}}^{2}\ \ \ \Bar{\xi}_{1}\frac{dx}{x}+\Bar{\xi}_{2}\frac{dy}{y}\mapsto (\Bar{D}_{x}(\Bar{\xi}_{2})-\Bar{D}_{y}(\Bar{\xi}_{1}))\frac{dx}{x}\wedge\frac{dy}{y}.     
\end{align*}
Note that $H_{x}=x\frac{\partial f_{t}}{\partial x}$, $H_{y}=y\frac{\partial f_{t}}{\partial y}$. $x\frac{\partial}{\partial x}(\Bar{S}^{i})\subseteq\Bar{S}^{i}$ and $y\frac{\partial}{\partial y}(\Bar{S}^{i})\subseteq\Bar{S}^{i}$. Due to the non-degeneracy of our family, we have the following theorem on the two cohomological spaces:
\begin{thm} \label{theorem 2.4}
\normalfont(Haessig and Sperber\cite{DS01}, theorem 2.2) \textit{For every $t\in\mathbb{F}_{q}^{*}$, both $(\Omega_{\Bar{S}}^{\bullet},\triangledown(H))$ and $(\Omega_{\Bar{S}}^{\bullet},\triangledown(\Bar{D}))$ are acyclic except in the top dimension $\mathrm{2}$. In both cases, $H^{2}$ is a finitely free $\mathbb{F}_{q}$-algebra of rank $2n+1$. For each $i\in\mathbb{Z}_{\geq 0}$ we choose a monomial basis $B_{i}$ consisting of monomials of weight $i/n$ for an $\mathbb{F}_{q}$-vector space $V_{i}$ such that the $i$-th graded piece $\Bar{S}^{i}$ of $\Bar{S}$ may be written as
$$\Bar{S}^{i}=V_{i}\oplus( H_{x}\Bar{S}^{i-n} +H_{y}\Bar{S}^{i-n}).$$
We write $B=\mathop{\cup}\limits_{i\geq 0}B_{i}$, if
 $V=\mathop{\sum}\limits_{i\geq 0}V_{i}$ is a $\mathbb{F}_{q}$-vector space with basis $B$, then we have
$$H^{2}(\Omega_{\Bar{S}}^{\bullet},\triangledown(H))=\Bar{S}/(H_{x}\Bar{S}+ H_{y}\Bar{S})\simeq V$$
as well that
$$H^{2}(\Omega_{\Bar{S}}^{\bullet},\triangledown(\Bar{D}))=\Bar{S}/(\Bar{D}_{x}\Bar{S}+ \Bar{D}_{y}\Bar{S})\simeq V.$$}
\end{thm}
\par We begin with a lemma which will be helpful in computing the cohomology and the Hodge polygon: 
\begin{lem} \label{lemma 2.5}
\textit{For every $(a,b)\in\mathbb{Z}^{2}$ with $\omega(a,b)=i/n$, we have:\\
(a) when $0\leq i\leq n-1$, $(a,b)=(i,0)$,\\
(b) when $n\leq i\leq 2n-1$, $(a,b)=(i-n,1),(i,0)$ or $(i-n-1,-1)$,\\
(c) when $i=2n$, $(a,b)=(n,1),(2n,0),(0,2),(-1,0),(-2,-2),(n-1,-1)$.}
\end{lem}
\textit{Proof:} Combinatorially, we can fit in all the $\mathbb{Z}$-lattice points in $\frac{i}{n}\triangle$ and find the number of lattice points on the boundary for each $i\in\mathbb{Z}_{\geq 0}$, then the lemma will be seen by the value of $W_{\triangle}(i)$, number of intersection points of $M(\triangle)$ and $\frac{i}{n}\triangle$. $\hfill\square$ \\
\par We set the notation 
\begin{equation}
\varepsilon_{i}=\left\{\begin{array}{cc} 
   x^{i}  &  \mathrm{when}\ 0\leq i\leq n ,\\
    x^{i-n}y &  \mathrm{when}\  n+1\leq i\leq 2n. \label{2.16}
\end{array}\right.    
\end{equation}
Denote $\varepsilon_{i}=x^{\varepsilon_{i}(x)}y^{\varepsilon_{i}(y)}$, we see $\omega(\varepsilon_{i}(x),\varepsilon_{i}(y))=i/n$. With an abuse using of notation, $\varepsilon_{i}$ also represents $(\varepsilon_{i}(x),\varepsilon_{i}(y))$ in all the following arguments.
\begin{thm} \label{theorem 2.6}
\textit{$\{\varepsilon_{i}\}_{0\leq i\leq 2n}$ is a basis for $H^{2}(\Omega_{\Bar{S}}^{\bullet},\triangledown(\Bar{D})).$ Precisely speaking, we have $$\Bar{S}=\mathop{\oplus}\limits_{i=0}^{2n}\mathbb{F}_{q}\varepsilon_{i}\oplus(\Bar{D}_{x}\Bar{S}+\Bar{D}_{y}\Bar{S}).$$
And moreover, for any $i\geq 0$, if $\mu\in\mathop{\oplus}\limits_{j=0}^{i}\Bar{S}^{(j)}=\mathrm{Fil}^{i}(\Bar{S})$, we have $$\mu=\sum_{j=0}^{\mathrm{min}\{i,2n\}}\Bar{a}(\mu,\varepsilon_{j})\varepsilon_{j}+\Bar{D}_{x}\Bar{\zeta}_{x}(\mu)+\Bar{D}_{y}\Bar{\zeta}_{y}(\mu)$$ for some $\Bar{a}(\mu,\varepsilon_{j})\in\mathbb{F}_{q}$, $\Bar{\zeta}_{x}(\mu),\ \Bar{\zeta}_{y}(\mu)\in\mathrm{Fil}^{i-n}(\Bar{S}).$} 
\end{thm}
\textit{Proof:} By Wang and Yang \cite{WY01} section 3, $\{\varepsilon_{i}\}_{0\leq i\leq 2n}$ is a basis for $H^{2}(\Omega^{\bullet}_{\Bar{S}},\triangledown(H))$, then theorem 2.3 shows it is also a basis for $H^{2}(\Omega_{\Bar{S}}^{\bullet},\triangledown(\Bar{D})).$  $\hfill\square$ \\
\par The next goal is to get a basis for $H^{2}(\Omega_{\mathcal{C}_{0}}^{\bullet}, \triangledown(D))$, we use the reduction map Pr defined in \hyperref[2.13]{(2.13)} as a bridge passing from the reduction cohomology on $\Bar{S}$ back to the cohomology on $\mathcal{C}_{0}$. For $(a,b)\in\mathbb{Z}^{2}$, we denote $(\widetilde{a,b}):=\widetilde{\gamma}^{n\omega(a,b)}x^{a}y^{b}$, this notational convention shows up throughout the subsequent material. Using this notation, we have the following result: 
\begin{thm} \label{theorem 2.7}
\textit{$\{\widetilde{\varepsilon}_{i}\}_{0\leq i\leq 2n}$ is a basis for $H^{2}(\Omega_{\mathcal{C}_{0}}^{\bullet}, \triangledown(D)).$ More precisely, we have $$\mathcal{C}_{0}=\mathop{\oplus}\limits_{i=0}^{2n}\mathbb{Z}_{q}[\widetilde{\gamma}]\widetilde{\varepsilon}_{i}\oplus(D_{x}\mathcal{C}_{0}+D_{y}\mathcal{C}_{0}).$$}
\end{thm}
\textit{Proof.} We just need to show for any $\eta\in\mathcal{C}_{0}$, there exists $\{a_{i}(\eta)\}_{0\leq i\leq 2n}\subseteq\mathbb{Z}_{q}[\widetilde{\gamma}]$ and $\xi_{x}(\eta)$, $\xi_{y}(\eta)\subseteq\mathcal{C}_{0}$ such that 
\begin{equation}
\eta=\sum_{i=0}^{2n}a_{i}(\eta)\widetilde{\varepsilon_{i}}+D_{x}(\xi_{x}(\eta))+D_{y}(\xi_{y}(\eta)).    
\end{equation}
\par Let $\Bar{\eta}=\mathrm{Pr}(\eta)$, then by theorem 2.5, we have the expression
\begin{equation}
 \Bar{\eta}=\sum_{i=0}^{2n}\Bar{\alpha}_{i}^{(1)}(\eta)\varepsilon_{i}+\Bar{D}_{x}\Bar{\xi}_{x}^{(1)}(\eta)+\Bar{D}_{y}\Bar{\xi}_{y}^{(1)}(\eta)   
\end{equation}
where $\Bar{\alpha}_{i}^{(1)}\in\mathbb{F}_{q}$ and $\Bar{\xi}_{x}^{(1)}(\eta),\ \Bar{\xi}_{y}^{(1)}(\eta)\in\Bar{S}$. Now choose $\alpha_{i}^{(1)}(\eta)$ as a Teichmüller for $\Bar{\alpha}_{i}^{(1)}(\eta)$ in $\mathbb{Z}_{q}$, and choose some $\xi_{x}^{(1)}(\eta),\ \xi_{y}^{(1)}(\eta)$ as the preimages of $\Bar{\xi}_{x}^{(1)}(\eta),\ \Bar{\xi}_{y}^{(1)}(\eta)$ in $\mathcal{C}_{0}$ via the reduction map Pr:
$$\Bar{\eta}=\mathrm{Pr}(\eta)=\mathrm{Pr}(\sum_{i=0}^{2n}\alpha_{i}^{(1)}\widetilde{\varepsilon}_{i}+D_{x}\xi_{x}^{(1)}(\eta)+D_{y}\xi_{y}^{(1)}(\eta)).$$
Since $\mathcal{C}_{0}/\widetilde{\gamma}\mathcal{C}_{0}\mathop{\simeq}\limits^{\mathrm{Pr}}\Bar{S}$, we have $\eta-(\sum\limits_{i=0}^{2n}\alpha_{i}^{(1)}\widetilde{\varepsilon}_{i}+D_{x}\xi_{x}^{(1)}(\eta)+D_{y}\xi_{y}^{(1)}(\eta))=\widetilde{\gamma}\eta^{(1)}$ for some $\eta^{(1)}\in\mathcal{C}_{0}$. Recursively applying above procedure we will get 
$$\eta^{(k-1)}-(\sum_{i=0}^{2n}\alpha_{i}^{(k)}\widetilde{\varepsilon}_{i}+D_{x}\xi_{x}^{(k)}(\eta)+D_{y}\xi_{y}^{(k)}(\eta))=\widetilde{\gamma}\eta^{(k)}$$
for some $\eta^{(k)}\in\mathcal{C}_{0}$. Let $a_{i}(\eta)=\sum\limits_{k\geq 0}\alpha_{i}^{(k)}(\eta)\widetilde{\gamma}^{k}$, $\xi_{x}(\eta)=\sum\limits_{k\geq 0}\xi_{x}^{(k)}(\eta)\widetilde{\gamma}^{k}$ and $\xi_{y}(\eta)=\sum\limits_{k\geq 0}\xi_{y}^{(k)}(\eta)\widetilde{\gamma}^{k}$. As $k\rightarrow\infty$, $|\widetilde{\gamma}^{k}|_{p}\rightarrow 0$, so the sum for $a_{i}(\eta)$ converges $\widetilde{\gamma}$-adically, then $a_{i}(\eta)\in\mathbb{Z}_{q}[\widetilde{\gamma}]$. Similar reason, $\xi_{x}(\eta)$ and $\xi_{y}(\eta)$ are well-defined under the superior norm of $\mathcal{C}_{0}$. By the recursive relations, they are the elements fitting into the equation (2.17). $\hfill\square$ 
\par Using $\{(\widetilde{a,b})\}_{(a,b)\in\mathbb{Z}^{2}}$ as an orthonormal basis for $\mathcal{C}_{0}$, let $A((\widetilde{c,d}),(\widetilde{a,b}))$ denote the coefficient of $(\widetilde{c,d})$ in the expression $\alpha_{1}((\widetilde{a,b}))=\sum_{(c,d)\in\mathbb{Z}^{2}}A((\widetilde{c,d}),(\widetilde{a,b}))\cdot(\widetilde{c,d})$
A simple calculation shows:
\begin{equation}
A((\widetilde{c,d}),(\widetilde{a,b}))=B^{\sigma^{-1}}(pc-a,pd-b)\widetilde{\gamma}^{n\omega(a,b)-n\omega(c,d)}    \label{2.19} 
\end{equation}
where $B^{\sigma^{-1}}$ means applying $\sigma^{-1}$ to the coefficients of $B$. 
\par Without any further conditions, we have the $p$-adic estimation as follows. 
\begin{lem} \label{lemma 2.8}
\textit{For any $(a,b)$, $(c,d)\in\mathbb{Z}^{2}$, $\ord_{p}A((\widetilde{c,d}),(\widetilde{a,b}))\geq\omega(c,d)$.}
\end{lem}
\textit{Proof.} By \hyperref[lemma 2.3]{lemma 2.3} and the triangle inequality for the weight function we easily obtain the result.
$\hfill\square$ \\
\par To get a better $p$-adic estimation, denote $0<\varpi<n$ the integer such that $p\varpi\equiv 1\ \mathrm{mod}\ n$. Let
\begin{equation}
g(i)=\left\{\begin{array}{cc}
i  &  \mathrm{when}\ i=0,n,2n, \\
\varpi i+n\lceil\frac{-\varpi i}{n}\rceil  & \ \ \ \mathrm{when}\ 1\leq i\leq n-1, \\
n+\varpi i+n\lceil\frac{-\varpi i}{n}\rceil  & \ \ \ \ \ \ \ \ \ \ \mathrm{when}\ n+1\leq i\leq 2n-1.    \label{2.20}
\end{array}\right.    
\end{equation}
Here $\lceil\ \rceil$ is the ceiling function. We use this notation in the following document.
\par Notice that $g$ is a bijection from $\mathbb{Z}\cap[0,2n]$ to itself. $0\leq g(i)\leq n-1$ when $0\leq i\leq n-1$, $n+1\leq g(i)\leq 2n-1$ when $n+1\leq i\leq 2n-1$, and more importantly, $g(i)\equiv \varpi i\ \mathrm{mod}\ n$ for $0\leq i\leq 2n$.
We then have the following theorem which refines lemma 2.7 to some extent. 
\begin{thm} \label{theorem 2.9}
\textit{For any $(c,d)\in\mathbb{Z}^{2}$, $0\leq i\leq 2n$, if $\omega(c,d)\leq \frac{g(i)}{n}$ and $(c,d)\neq\varepsilon_{g(i)}$, then $\ord_{p}A((\widetilde{c,d}),\widetilde
{\varepsilon_{i}})>\omega(c,d)$. } 
\end{thm}
\textit{Proof.} Firstly, by \hyperref[2.19]{(2.19)} and \hyperref[lemma 2.3]{lemma 2.3} we have
\begin{align*}
\ord_{p}A((\widetilde{c,d}),\widetilde
{\varepsilon_{i}})&=\frac{\frac{i}{n}-\omega(c,d)}{p-1}+\ord_{p}B^{\sigma^{-1}}(pc-\varepsilon_{i}(x),pd-\varepsilon_{i}(y))\\
&\geq\frac{\frac{i}{n}-\omega(c,d)}{p-1}+\frac{k+l+m}{p-1}\\
&=\frac{\frac{i}{n}-\omega(c,d)}{p-1}+\frac{\frac{2n+1}{n}k+\frac{pc}{n}+pd-\frac{i}{n}}{p-1}\\
&=\omega(c,d)+\frac{2n+1}{n}(k-pm(c,d))
\end{align*}
where $k$ is the smallest non-negative integer such that the triple $(k,l,m)$ lies in $$I(pc-\varepsilon_{i}(x),pd-\varepsilon_{i}(y))=\{(k,l,m)\in\mathbb{Z}_{\geq 0}^{3}|nl-k=pc-\varepsilon_{i}(x),\ m-k=pd-\varepsilon_{i}(y)\}$$
as we defined in \hyperref[lemma 2.3]{lemma 2.3}.
\par We aim to prove the theorem by showing that for this smallest $k$, $k>pm(c,d)$. \\
\textit{Case I.} $m(c,d)=0$. 
\par In this case, $\omega(c,d)=\frac{c}{n}+d\leq\omega(\varepsilon_{g(i)})=\frac{g(i)}{n}$, $c\geq 0$ and $d\geq 0$. We show by controdiction via setting $k=0$. Then $nl=pc-\varepsilon_{i}(x)$, $m=pd-\varepsilon_{i}(y)$. 
\par If $c=0$, $nl=-\varepsilon_{i}(x)$ for $l\in\mathbb{Z}_{\geq 0}$, we must have $\varepsilon_{i}(x)=0$, the only one $\varepsilon_{i}$ with $\varepsilon_{i}(x)=0$ is $\varepsilon_{0}$, so $i=0$, $g(i)=0$. Then $\omega(c,d)\leq 0$, the only choice is $(c,d)=(0,0)$, which violates the condition $(c,d)\neq\varepsilon_{0}$. So we must have $c\geq 1$.
\par If $d\geq 2$, we have $\frac{g(i)}{n}\geq\frac{c}{n}+d\geq 2+\frac{1}{n}$, $g(i)\geq 2n+1$. This controdicts with the definition of $g(i)$. Therefore $d\leq 1$. $m=pd-\varepsilon_{i}(y)$ with $m\geq 0$ and $\varepsilon_{i}(y)=0\ \mathrm{or}\ 1$, this means $d$ cannot be negative, so we must have $d=0$ or 1.\\
\textit{Subcase I.1.} $d=0$.
\par In this subcase, $\omega(c,0)=\frac{c}{n}\leq\frac{g(i)}{n}$, so $c\leq g(i)$. $m=pd-\varepsilon_{i}(y)=-\varepsilon_{i}(y)\geq 0$, so $\varepsilon_{i}(y)=0$, this means $0\leq i\leq n$, and $0\leq g(i)\leq n$.
\par Notice that $g(i)\equiv \varpi i\ \mathrm{mod}\ n$ where $p\varpi\equiv 1\ \mathrm{mod}\ n$, we have 
\begin{equation}
    pg(i)\equiv i\ \mathrm{mod}\ n\ \ \ \mathrm{for}\ \mathrm{all}\ 0\leq i\leq 2n.
\end{equation}
\par When $0\leq i\leq n$, $\varepsilon_{i}=x^{i}$, so $\varepsilon_{i}(x)=i$, $nl=pc-i$. $p\nmid n$, $nl=pc-i=p(c-g(i))+pg(i)-i$ implies $g(i)-c\equiv 0\ \mathrm{mod}\ n$. But $0\leq g(i)-c<g(i)\leq n$, the only choice is $g(i)=c$. Therefore we have $(c,d)=(g(i),0)=\varepsilon_{g(i)}$, violating the condition $(c,d)\neq\varepsilon_{g(i)}$.\\
\textit{Subcase I.2.} $d=1$.
\par In this subcase, $\frac{1}{n}+1\leq\frac{c}{n}+1=\omega(c,1)\leq\frac{g(i)}{n}$, so $n+1\leq g(i)$ and $c\leq g(i)-n$. Moreover we have $1+n\leq i,\  g(i)\leq 2n$. When $n+1\leq i\leq 2n$, we see $\varepsilon_{i}=x^{i-n}y$, so $\varepsilon_{i}(x)=i-n$, $\varepsilon_{i}(y)=1$. We also have $nl=pc-i+n=p(c-g(i))+pg(i)-i+n$, which implies $g(i)-c\equiv 0\ \mathrm{mod}\ n$. But $n\leq g(i)-c<g(i)\leq 2n$, the only choice is $g(i)-c=n$. Therefore we also have $(c,d)=(g(i)-n,1)=\varepsilon_{g(i)}$, again violating $(c,d)\neq\varepsilon_{g(i)}$.\\
\textit{Case II.} $m(c,d)=-c$.
\par In this case, $c\leq 0$ and $c\leq d$. We need to show $k>-pc$. Since $k=nl-pc+\varepsilon_{i}(x)\geq -pc$, the smallest $k$ could be $-pc$. Suppose $k=-pc$, then $nl-k=pc-\varepsilon_{i}(x)$ implies $nl=-\varepsilon_{i}(x)\leq 0$. The only choice is $i=0$, then $g(i)=0$, and $\omega(c,d)\leq\omega(\varepsilon_{0})=0$, which gives us $(c,d)=(0,0)$, violating the condition $(c,d)\neq\varepsilon_{0}$.\\
\textit{Case III.} $m(c,d)=-d$.
\par In this case $d\leq 0$ and $d\leq c$. We need to show $k>-pd$. Since $k=m-pd+\varepsilon_{i}(y)\geq-pd$, the smallest $k$ could be $-pd$. Again suppose $k=-pd$, then $m-k=pd-\varepsilon_{i}(y)$ implies $m=-\varepsilon_{i}(y)\geq 0$. So $\varepsilon_{i}(y)=0$, we have $0\leq i,\ g(i)\leq n$ and $\varepsilon_{i}(x)=n$. $nl=pc-pd-i=p(c-d-g(i))+pg(i)-i$ implies $g(i)+d-c\equiv 0\ \mathrm{mod}\ n$. But $\omega(c,d)=\frac{c}{n}+d-\frac{2n+1}{n}d\leq\frac{g(i)}{n}$ gives $c-d-nd\leq g(i)$, which implies $0\leq -nd\leq g(i)+d-c\leq g(i)\leq n$. So we only have two cases.\\
\textit{Subcase III.1.} $g(i)+d-c=0$.
\par By $0\leq -nd\leq g(i)+d-c=0$, we see $d=0$. This goes back to \textit{Case I.} and the controdiction follows.\\
\textit{Subcase III.2.} $g(i)+d-c=n$.
\par By $g(i)+d-c\leq g(i)\leq n$ we have $d=c$ and $g(i)=n$. Therefore $i=n$, so $nl=pc+k-n=pc-pd-n=-n$, $l=-1$, which violates $l\geq 0$.
\par We see in all cases $k>pm(c,d)$ is true, therefore the theorem follows. $\hfill\square$ \\
\par We then have the following $p$-adic estimations of $A(\widetilde{\varepsilon}_{j},\widetilde{\varepsilon}_{i})$:
\begin{thm} \label{theorem 2.10}
\textit{Suppose $p>n$. For $0\leq i\leq 2n$, we have $\ord_{p}A(\widetilde{\varepsilon}_{g(i)},\widetilde{\varepsilon}_{i})=\frac{g(i)}{n}$. For $0\leq i\leq n$, $0\leq j\leq n-1$, we have $\ord_{p}A(\widetilde{\varepsilon}_{j},\widetilde{\varepsilon}_{i})=\frac{j}{n}+\frac{(2n+1)\alpha_{i,j}}{n(p-1)}$.}
\end{thm}
\textit{Proof.} For the first statement,
$$\ord_{p}A(\widetilde{\varepsilon}_{g(i)},\widetilde{\varepsilon}_{i})=\frac{i-g(i)}{n(p-1)}+\ord_{p}B^{\sigma^{-1}}(p\varepsilon_{g(i)}(x)-\varepsilon_{i}(x),p\varepsilon_{g(i)}(y)-\varepsilon_{i}(y)),$$
and $$\ord_{p}B^{\sigma^{-1}}(p\varepsilon_{g(i)}(x)-\varepsilon_{i}(x),p\varepsilon_{g(i)}(y)-\varepsilon_{i}(y))\geq\mathrm{inf}\frac{k+l+m}{p-1}$$
where $(k,l,m)\in I(p\varepsilon_{g(i)}(x)-\varepsilon_{i}(x),p\varepsilon_{g(i)}(y)-\varepsilon_{i}(y))$ as defined in \hyperref[lemma 2.3]{lemma 2.3}. So we have
$$nl-k=p\varepsilon_{g(i)}(x)-\varepsilon_{i}(x), m-k=p\varepsilon_{g(i)}(y)-\varepsilon_{i}(y).$$
The inequality is an equality if and only if we can find only one triple $(k.l.m)$ such that $\ord_{p}(a_{k}a_{l}a_{m})$ is strictly the smallest one. Here $a_{k},a_{l},a_{m}$ are the coefficients of the splitting function $\Theta_{\infty}$. Recall $\ord_{p}a_{j}=\frac{j}{p-1}$ if $0\leq j\leq p-1$, and in general $\ord_{p}a_{j}\geq\frac{j}{p-1}$.
\par When $0\leq i\leq n$, $0\leq g(i)\leq n$, so $\varepsilon_{i}=x^{i}$ and $\varepsilon_{g(i)}=x^{g(i)}$. Therefore we have 
$$nl-k=pg(i)-i,\ m-k=0.$$
Recall $n|pg(i)-i$, we see $k=m=0,\ l=\frac{pg(i)-i}{n}$ is the case that $k$ is the smallest. Since $p>n$, we see $0\leq \frac{pg(i)-i}{n}\leq p-1$, so we have the accurate $p$-adic estimation
$$\ord_{p}(a_{0}^{2}a_{\frac{pg(i)-i}{n}})=\frac{pg(i)-i}{n(p-1)}.$$
For all other triples $(k,l,m)\in I(p\varepsilon_{g(i)}(x)-\varepsilon_{i}(x),p\varepsilon_{g(i)}(y)-\varepsilon_{i}(y))$, $\mathrm{Ord}_{p}(a_{k}a_{l}a_{m})$ is strictly larger than $\frac{pg(i)-i}{n(p-1)}.$ So we get the accurate estimation
$$\ord_{p}A(\widetilde{\varepsilon}_{g(i)},\widetilde{\varepsilon}_{i})=\frac{i-g(i)}{n(p-1)}+\frac{pg(i)-i}{n(p-1)}=\frac{g(i)}{n}.$$
\par When $n+1\leq i\leq 2n$, $\varepsilon_{i}=x^{i-n}y$ and $\varepsilon_{g(i)}=x^{g(i)-n}y$, in this case we have
$$nl-k=pg(i)-i+n-pn,\ m-k=p-1.$$
Again to make $k$ the smallest, we set $k=0$, then $m=p-1$, $l=\frac{pg(i)-i}{n}+1-p$. By $p>n$ we also have $0\leq\frac{pg(i)-i}{n}+1-p\leq p-1$. Then same as the case when $0\leq i\leq n$, the smallest $p$-adic order is
$$\ord_{p}(a_{0}a_{p-1}a_{\frac{pg(i)-i}{n}+1-p})=\frac{pg(i)-i}{n(p-1)}$$ and the first statement follows.
\par For the second statement, $0\leq i\leq n$, $0\leq j\leq n-1$, we have
$$\ord_{p}A(\widetilde{\varepsilon}_{j},\widetilde{\varepsilon}_{i})=\frac{i-j}{n(p-1)}+\ord_{p}B^{\sigma^{-1}}(pj-i,0)$$
and $\ord_{p}B^{\sigma^{-1}}(pj-i,0)\geq\mathop{\mathrm{inf}}\limits_{(k,l,m)\in I(pj-i,0)}\ord_{p}(a_{k}a_{l}a_{m})$. The triples $(k,l,m)$ satisfies
$$nl-k=pj-i,\ m-k=0.$$
Then the choice of the smallest $k$ is $k=\alpha_{i,j}$ as defined above. Again $p>n$ shows that $\alpha_{i,j}\leq p-1$. And $0\leq l=\frac{\alpha_{i,j}+pj-i}{n}=\lceil\frac{pj-i}{n}\rceil\leq p-1$ when $0\leq j\leq n-1$. The smallest $p$-adic order is
\begin{align*}
\ord_{p}(a_{\alpha_{i,j}}^{2}a_{\frac{\alpha_{i,j}+pj-i}{n}})=\frac{2\alpha_{i,j}}{p-1}+\frac{\alpha_{i,j}+pj-i}{n(p-1)}  . 
\end{align*}
Plug in the expression of $\ord_{p}A(\widetilde{\varepsilon}_{j},\widetilde{\varepsilon}_{i})$ the second statement follows.  $\hfill\square$
\par The next goal is to establish the explicit formulas for the coefficients of the Frobenius action in the cohomological level. By \hyperref[theorem 2.7]{theorem 2.7} $\{\widetilde{\varepsilon}_{i}\}_{0\leq i\leq 2n}$ forms a basis of $H^{2}(\Omega_{\mathcal{C}_{0}}^{\bullet})=\mathcal{C}_{0}/(D_{x}\mathcal{C}_{0}+D_{y}\mathcal{C}_{0})$. Suppose $\widetilde{A}(\widetilde{\varepsilon}_{j},\widetilde{\varepsilon}_{i})$ is the coefficient of $\widetilde{\varepsilon}_{j}$ when we express $\alpha_{1}(\widetilde{\varepsilon}_{i})$ using the basis $\{\widetilde{\varepsilon}_{i}\}_{0\leq i\leq 2n}$, which means
\begin{equation}
\alpha_{1}(\widetilde{\varepsilon}_{i})=\sum_{j=0}^{2n}\widetilde{A}(\widetilde{\varepsilon}_{j},\widetilde{\varepsilon}_{i})\widetilde{\varepsilon}_{j}\ \mathrm{in}\ H^{2}(\Omega_{\mathcal{C}_{0}}^{\bullet}).  \label{2.23}
\end{equation}
\par Using this basis, for any $\beta\in\mathcal{C}_{0}$, by \hyperref[theorem 2.7]{theorem 2.7} we may write
\begin{equation}
\beta=\sum_{k=0}^{2n}a(\beta,\widetilde{\varepsilon}_{k})\widetilde{\varepsilon}_{k}+D_{x}\zeta_{x}(\beta)+D_{y}\zeta_{y}(\beta)    \label{2.24}
\end{equation}
where $a(\beta,\widetilde{\varepsilon}_{k})\in\mathbb{Z}_{q}[\widetilde{\gamma}]$ is unique, and $\zeta_{x}(\beta),\zeta_{y}(\beta)\in\mathcal{C}_{0}$. Therefore we have the following expression for $\alpha_{1}(\widetilde{\varepsilon}_{i})$ in $\mathcal{C}_{0}$:
\begin{align*}
\alpha_{1}(\widetilde{\varepsilon}_{i})&=\sum_{(c,d)\in\mathbb{Z}^{2}}A((\widetilde{c,d}),\widetilde{\varepsilon}_{i})(\widetilde{c,d})\\
&=\sum_{(c,d)\in\mathbb{Z}^{2}}\{A((\widetilde{c,d}),\widetilde{\varepsilon}_{i})\cdot[\sum_{j=0}^{2n}a((\widetilde{c,d}),\widetilde{\varepsilon}_{j})\widetilde{\varepsilon}_{j}+D_{x}\zeta_{x}((\widetilde{c,d}))+D_{y}\zeta_{y}((\widetilde{c,d})) ]\}  \\
&=\sum_{j=0}^{2n}[\sum_{(c,d)\in\mathbb{Z}^{2}}A((\widetilde{c,d}),\widetilde{\varepsilon}_{i})a((\widetilde{c,d}),\widetilde{\varepsilon}_{j})]\widetilde{\varepsilon}_{j}+D_{x}[\sum_{(c,d)\in\mathbb{Z}^{2}}A((\widetilde{c,d}),\widetilde{\varepsilon}_{i})\zeta_{x}((\widetilde{c,d}))] \\ &+D_{y}[\sum_{(c,d)\in\mathbb{Z}^{2}}A((\widetilde{c,d}),\widetilde{\varepsilon}_{i})\zeta_{y}((\widetilde{c,d}))]    .
\end{align*}
Compare this with \hyperref[2.23]{(2.23)}, we obtain the expression on the cohomological level:
\begin{equation}
\widetilde{A}(\widetilde{\varepsilon}_{j},\widetilde{\varepsilon}_{i})=\sum_{(c,d)\in\mathbb{Z}^{2}}A((\widetilde{c,d}),\widetilde{\varepsilon}_{i})a((\widetilde{c,d}),\widetilde{\varepsilon}_{j}).  \label{2.25}  
\end{equation}
Previous theorems give much estimations on $A((\widetilde{c,d}),\widetilde{\varepsilon}_{i})$. To study the Frobenius coefficients, we need to give some $p$-adic estimations on $a((\widetilde{c,d}),\widetilde{\varepsilon}_{j})$. Keep the convention that $\varepsilon_{j}=0$ for $j>2n$. Inspired by \cite{AS02} we proceed this by the following lemmas.
\begin{lem} \label{lemma 2.11}
\textit{Let $T^{(i)}$ be the $\mathbb{Z}_{q}[\widetilde{\gamma}]$-submodule of $\mathcal{C}_{0}$ generated by $\{(\widetilde{c,d})\}_{(c,d)\in\mathbb{Z}^{2},n\omega(c,d)\leq i}$. Let $D_{x}^{(1)}=x\frac{\partial}{\partial x}+\gamma(nx^{n}-\frac{\Bar{t}}{xy})$, $D_{y}^{(1)}=y\frac{\partial}{\partial y}+\gamma(y-\frac{\Bar{t}}{xy})$ where $\Bar{t}$ is a Teichmüller of $t$.\\
(a) If $0\leq i\leq n-1$, then $T^{(i)}$ is generated by $\{\widetilde{\varepsilon}_{j}\}_{0\leq j\leq i}$. \\
(b) If $n\leq i$, then for any $\beta\in T^{(i)}$, there exist $\{a'(\beta,\widetilde{\varepsilon}_{j})\}_{0\leq j\leq i}\subseteq\mathbb{Z}_{q}[\widetilde{\gamma}]$, $\zeta'_{x}(\beta),\zeta'_{y}(\beta)\in T^{(i-n)}$ such that $$\beta=\sum_{j=0}^{i}a'(\beta,\widetilde{\varepsilon}_{j})\widetilde{\varepsilon}_{j}+D_{x}^{(1)}\zeta'_{x}(\beta)+D_{y}^{(1)}\zeta'_{y}(\beta).$$} 
\end{lem}
\textit{Proof.} The first statement follows from \hyperref[lemma 2.5]{lemma 2.5} immediately. For the second statement, note that under the projection Pr, the reduction of $D_{x}$ (resp. $D_{y}$) and $D^{(1)}_{x}$ (resp. $D^{(1)}_{y}$) are the same operator, $\Bar{D}_{x}$ (resp. $\Bar{D}_{y}$) as defined in \hyperref[2.14]{(2.14)} (resp. \hyperref[2.15]{(2.15)}). Similar like $\Bar{D}_{x}$ and $\Bar{D}_{y}$ we have $D^{(1)}_{x}(T^{(i-n)})\subseteq T^{(i)}$ and $D^{(1)}_{y}(T^{(i-n)})\subseteq T^{(i)}$.
\par Reduce $\beta$ modulo $\widetilde{\gamma}$, we have $\Bar{\beta}=\mathrm{Pr}(\beta)\in\mathop{\oplus}_{j=0}^{i}\Bar{S}^{(j)}$, then by \hyperref[theorem 2.6]{theorem 2.6}, 
$$\Bar{\beta}=\sum_{j=0}^{i}\Bar{a}'^{(1)}(\beta,\widetilde{\varepsilon}_{j})\widetilde{\varepsilon}_{j}+D_{x}^{(1)}\Bar{\zeta}'^{(1)}_{x}(\beta)+D_{y}^{(1)}\Bar{\zeta}'^{(1)}_{y}(\beta)$$
for some $\{\Bar{a}'^{(1)}(\beta,\widetilde{\varepsilon}_{j})\}_{0\leq j\leq i}\subseteq\mathbb{F}_{q}$, $\Bar{\zeta}'^{(1)}_{x}(\beta),\Bar{\zeta}'^{(1)}_{y}(\beta)\in \mathop{\oplus}\limits_{j=0}^{i-n}\Bar{S}^{(j)}$. Lifting this back to $\mathcal{C}_{0}$, using above commutativities we obtain
$$\beta=\sum_{j=0}^{i}a'^{(1)}(\beta,\widetilde{\varepsilon}_{j})\widetilde{\varepsilon}_{j}+D_{x}^{(1)}\zeta'^{(1)}_{x}(\beta)+D_{y}^{(1)}\zeta'^{(1)}_{y}(\beta)+\widetilde{\gamma}\beta^{(1)}$$
with some $\{a'^{(1)}(\beta,\widetilde{\varepsilon}_{j})\}_{0\leq j\leq i}\subseteq\mathbb{Z}_{q}[\widetilde{\gamma}]$, $\zeta'^{(1)}_{x}(\beta),\zeta'^{(1)}_{y}(\beta)\in T^{(i-n)}$ and $\beta^{(1)}\in T^{(i)}$. Repeat above procedure for $\beta^{(1)}$ and follow with exactly the same recursive argument as in \hyperref[theorem 2.7]{theorem 2.7} we obtain the lemma. $\hfill\square$ \\
\par If we set $T^{(0)}\simeq\mathbb{Z}_{q}[\widetilde{\gamma}]$ and $T^{(i)}=0$ for $i<0$, we see the two statements in above lemma actually state the same result. We then pass from $D_{x}^{(1)}$ (resp. $D_{y}^{(1)}$) to $D_{x}$ (resp. $D_{y}$) and keep this convention in the following lemmas:
\begin{lem} \label{lemma 2.12}
\textit{Suppose $p>n$. If $n\leq i$, then for any $\beta\in T^{(i)}$, there exist $\varrho_{j}(\beta)\in T^{(j)}$ for all $j\geq i+1$ such that \hyperref[lemma 2.11]{lemma 2.11 (b)} can be rewritten as $$\beta=\sum_{j=0}^{i}a'(\beta,\widetilde{\varepsilon}_{j})\widetilde{\varepsilon}_{j}+D_{x}\zeta'_{x}(\beta)+D_{y}\zeta'_{y}(\beta)+\sum_{j=i+1}^{\infty}p^{j+n-i-1}\varrho_{j}(\beta).$$} 
\end{lem}
\textit{Proof.} Recall $$D_{x}=D_{x}^{(1)}+\sum_{m=1}^{\infty}\gamma_{m}p^{m}(nx^{np^{m}}-\frac{\sigma^{m}(\Bar{t})}{x^{p^{m}}y^{p^{m}}}), \  D_{y}=D_{y}^{(1)}+\sum_{m=1}^{\infty}\gamma_{m}p^{m}(y^{p^{m}}-\frac{\sigma^{m}(\Bar{t})}{x^{p^{m}}y^{p^{m}}}).$$Then by \hyperlink{Lemma 2.10.}{lemma 2.10} we may write 
\begin{align*}
\beta&=\sum_{j=0}^{i}a'(\beta,\widetilde{\varepsilon}_{j})\widetilde{\varepsilon}_{j}+D_{x}^{(1)}\zeta'_{x}(\beta)+D_{y}^{(1)}\zeta'_{y}(\beta) \\
&=\sum_{j=0}^{i}a'(\beta,\widetilde{\varepsilon}_{j})\widetilde{\varepsilon}_{j}+D_{x}\zeta'_{x}(\beta)+D_{y}\zeta'_{y}(\beta)-\sum_{m=1}^{\infty}\gamma_{m}p^{m}(nx^{np^{m}}-\frac{\sigma^{m}(\Bar{t})}{x^{p^{m}}y^{p^{m}}})\zeta'_{x}(\beta)\\
&-\sum_{m=1}^{\infty}\gamma_{m}p^{m}(y^{p^{m}}-\frac{\sigma^{m}(\Bar{t})}{x^{p^{m}}y^{p^{m}}})\zeta'_{y}(\beta)\\
&=\sum_{j=0}^{i}a'(\beta,\widetilde{\varepsilon}_{j})\widetilde{\varepsilon}_{j}+D_{x}\zeta'_{x}(\beta)+D_{y}\zeta'_{y}(\beta)\\
&+\sum_{m=1}^{\infty}-\frac{\gamma_{m}p^{m}}{\gamma^{p^{m}}}[\gamma^{p^{m}}(nx^{np^{m}}-\frac{\sigma^{m}(\Bar{t})}{x^{p^{m}}y^{p^{m}}})\zeta'_{x}(\beta)+\gamma^{p^{m}}(y^{p^{m}}-\frac{\sigma^{m}(\Bar{t})}{x^{p^{m}}y^{p^{m}}})\zeta'_{y}(\beta)].
\end{align*}
\par Clearly we see $nx^{np^{m}}-\frac{\sigma^{m}(\Bar{t})}{x^{p^{m}}y^{p^{m}}},\ y^{p^{m}}-\frac{\sigma^{m}(\Bar{t})}{x^{p^{m}}y^{p^{m}}}\in T^{(p^{m})}$. Together with $\zeta'_{x}(\beta),\ \zeta'_{y}(\beta)\in T^{(i-n)}$ we obtain
$$\gamma^{p^{m}}(nx^{np^{m}}-\frac{\sigma^{m}(\Bar{t})}{x^{p^{m}}y^{p^{m}}})\zeta'_{x}(\beta)+\gamma^{p^{m}}(y^{p^{m}}-\frac{\sigma^{m}(\Bar{t})}{x^{p^{m}}y^{p^{m}}})\zeta'_{y}(\beta)\in T^{(p^{m}+i-n)}$$
for any $m\geq 1$. Also note that 
$$\ord_{p}(\frac{\gamma_{m}p^{m}}{\gamma^{p^{m}}})=p^{m}-1$$
for every $m\geq 1$. We may write
$$\sum_{m=1}^{\infty}-\frac{\gamma_{m}p^{m}}{\gamma^{p^{m}}}[\gamma^{p^{m}}(nx^{np^{m}}-\frac{\sigma^{m}(\Bar{t})}{x^{p^{m}}y^{p^{m}}})\zeta'_{x}(\beta)+\gamma^{p^{m}}(y^{p^{m}}-\frac{\sigma^{m}(\Bar{t})}{x^{p^{m}}y^{p^{m}}})\zeta'_{y}(\beta)]=\sum_{m=1}^{\infty}p^{p^{m}-1}\varrho_{p^{m}+i-n}(\beta)$$
for some $\varrho_{j}(\beta)\in T^{(j)}$. Let $j=p^{m}+i-n$, then $p^{m}-1=j+n-i-1$. And we see $j\geq i+1$ since $p>n$. We rewrite the sum as $\sum\limits_{j=i+1}^{\infty}p^{j+n-i-1}\varrho_{j}(\beta)$ then the lemma follows. $\hfill\square$ 
\par Next we interact the weights and $p$-adic filtrations and get the following key lemma:
\begin{lem} \label{lemma 2.13}
\textit{Suppose $p>n$. If $n\leq i$, then for any $\beta\in T^{(i)}$, there exist $\{\widetilde{a}(\beta,\widetilde{\varepsilon}_{j})\}_{j\geq i+1}\subseteq\mathbb{Z}_{q}[\widetilde{\gamma}],\ \zeta_{x}(\beta),\zeta_{y}(\beta)\in\mathcal{C}_{0}$ such that $$\beta=\sum_{j=0}^{i}a(\beta,\widetilde{\varepsilon}_{j})\widetilde{\varepsilon}_{j}+\sum_{j=i+1}^{\infty}p^{j+n-i-1}\widetilde{a}(\beta,\widetilde{\varepsilon}_{j})\widetilde{\varepsilon}_{j}+D_{x}\zeta_{x}(\beta)+D_{y}\zeta_{y}(\beta)$$ where $a(\beta,\widetilde{\varepsilon}_{j})$'s are the coefficients expressed in \hyperref[2.24]{(2.24)}.} 
\end{lem}
\textit{Proof.} Keep in mind the convention $\varepsilon_{j}=0$ for $j>2n$. For any non-negative integer $N$, we claim 
\begin{equation}
\begin{aligned}
\beta=&\sum_{j=0}^{i}a^{(N)}(\beta,\widetilde{\varepsilon}_{j})\widetilde{\varepsilon}_{j}+\sum_{j=i+1}^{N+i}p^{j+n-i-1}\widetilde{a}^{(N)}(\beta,\widetilde{\varepsilon}_{j})\widetilde{\varepsilon}_{j}\\ &+D_{x}\zeta_{x}^{(N)}(\beta)+D_{y}\zeta_{y}^{(N)}(\beta)+\sum_{j=N+i+1}^{\infty}p^{j+n-i-1}\varrho_{j}^{(N)}(\beta)   \label{2.25}  
\end{aligned}    
\end{equation}
for some $\{a^{(N)}(\beta,\widetilde{\varepsilon}_{j})\}_{0\leq j\leq i}\cup\{\widetilde{a}^{(N)}(\beta,\widetilde{\varepsilon}_{j})\}_{i+1\leq j\leq N+i}\subseteq\mathbb{Z}_{q}[\widetilde{\gamma}]$, $\zeta_{x}^{(N)}(\beta),\zeta_{y}^{(N)}(\beta)\in\mathcal{C}_{0}$, and $\varrho_{j}^{(N)}(\beta)\in T^{(j)}$ for $j\geq N+i+1$.
\par We now show this claim by induction. Clearly \hyperref[lemma 2.12]{lemma 2.12} begins the induction for $N=0$ by setting $a^{(0)}(\beta,\widetilde{\varepsilon}_{j})=a'(\beta,\widetilde{\varepsilon}_{j})$ for $0\leq j\leq i$, $\zeta_{x}^{(0)}(\beta)=\zeta'_{x}(\beta)$, $\zeta_{y}^{(0)}(\beta)=\zeta'_{y}(\beta)$ and $\varrho_{j}^{(0)}(\beta)=\varrho_{j}(\beta)$ for $j\geq i+1$. Suppose the claim holds for $N$, we apply \hyperref[lemma 2.12]{lemma 2.12} again for the term $\varrho_{N+i+1}^{(N)}(\beta)\in T^{(N+i+1)}$ and obtain
\begin{equation}
\begin{aligned}
\varrho_{N+i+1}^{(N)}(\beta)=&\sum_{j=0}^{N+i+1}a'(\varrho_{N+i+1}^{(N)}(\beta),\widetilde{\varepsilon}_{j})\widetilde{\varepsilon}_{j}+D_{x}\zeta'_{x}(\varrho_{N+i+1}^{(N)}(\beta)) \\
&+D_{y}\zeta'_{y}(\varrho_{N+i+1}^{(N)}(\beta))+\sum_{j=N+i+2}^{\infty}p^{j+n-N-i-2}\varrho_{j}(\varrho_{N+i+1}^{(N)}(\beta))  \label{2.26}
\end{aligned}    
\end{equation}
for some $\{a'(\varrho_{N+i+1}^{(N)}(\beta),\widetilde{\varepsilon}_{j})\}_{0\leq j\leq N+i+1}\subseteq\mathbb{Z}_{q}[\widetilde{\gamma}]$, $\zeta'_{x}(\varrho_{N+i+1}^{(N)}(\beta)),\zeta'_{y}(\varrho_{N+i+1}^{(N)}(\beta))\in\mathcal{C}_{0}$ and  \\ $\varrho_{j}(\varrho_{N+i+1}^{(N)}(\beta))\in T^{(j)}$ for all $j\geq N+i+2$. Then by substituting \hyperref[2.26]{(2.26)} back into \hyperref[2.25]{(2.25)} we show the claim for $N+1$ with
\begin{align*}
&a^{(N+1)}(\beta,\widetilde{\varepsilon}_{j})=a^{(N)}(\beta,\widetilde{\varepsilon}_{j})+p^{N+n}a'(\varrho_{N+i+1}^{(N)}(\beta),\widetilde{\varepsilon}_{j}) \ \mathrm{for}\ 0\leq j\leq i, \\
&\widetilde{a}^{(N+1)}(\beta,\widetilde{\varepsilon}_{j})=\widetilde{a}^{(N)}(\beta,\widetilde{\varepsilon}_{j})+p^{N-j+i+1}a'(\varrho_{N+i+1}^{(N)}(\beta),\widetilde{\varepsilon}_{j})\ \mathrm{for}\ i+1\leq j\leq N+i, \\
&\widetilde{a}^{(N+1)}(\beta,\widetilde{\varepsilon}_{N+i+1})=a'(\varrho_{N+i+1}^{(N)}(\beta),\widetilde{\varepsilon}_{N+i+1}),  \\
&\zeta_{x}^{(N+1)}(\beta)=\zeta_{x}^{(N)}(\beta)+p^{N+n}\zeta'_{x}(\varrho_{N+i+1}^{(N)}(\beta)),\ \zeta_{y}^{(N+1)}(\beta)=\zeta_{y}^{(N)}(\beta)+p^{N+n}\zeta'_{y}(\varrho_{N+i+1}^{(N)}(\beta)), \\
&\varrho_{j}^{(N+1)}(\beta)=\varrho_{j}^{(N)}(\beta)+p^{n-1}\varrho_{j}(\varrho_{N+i+1}^{(N)}(\beta))\ \mathrm{for}\ j\geq N+i+2.
\end{align*}
\par Same convergent arguments like in \hyperref[theorem 2.7]{theorem 2.7} we may take limits as $N\rightarrow\infty$, by the uniqueness for coefficients in \hyperref[2.24]{(2.24)}, we have
$$a(\beta,\widetilde{\varepsilon}_{j})=\mathop{\mathrm{lim}}\limits_{N\rightarrow\infty}a^{(N)}(\beta,\widetilde{\varepsilon}_{j})\ \mathrm{for}\ 0\leq j\leq i.$$
Then we let $$\widetilde{a}(\beta,\widetilde{\varepsilon}_{j})=\mathop{\mathrm{lim}}\limits_{N\rightarrow\infty}\widetilde{a}^{(N)}(\beta,\widetilde{\varepsilon}_{j})\ \mathrm{for}\ j\geq i+1,$$
$$\zeta_{x}(\beta)=\mathop{\mathrm{lim}}\limits_{N\rightarrow\infty}\zeta_{x}^{(N)}(\beta),\ \zeta_{y}(\beta)=\mathop{\mathrm{lim}}\limits_{N\rightarrow\infty}\zeta_{y}^{(N)}(\beta).$$
Also note $\sum\limits_{j\geq N+i+1}p^{j+n-i-1}\varrho_{j}^{(N)}(\beta)$ vanishes as $N\rightarrow\infty$, therefore the lemma follows. $\hfill\square$
\par With the above lemmas, we will have the following key result which gives a $p$-adic estimation for the coefficient in \hyperref[2.24]{(2.24)}:
\begin{lem} \label{lemma 2.14}
    \textit{ \\ (a) If $0\leq i\leq n-1$, for any $(c,d)\in\mathbb{Z}^{2}$ with $n\omega(c,d)\leq n-1$ we have
$$a((\widetilde{c,d}),\widetilde{\varepsilon}_{i})=\left\{\begin{array}{cc}
   1  &  \ \ \ \  \ \  \mathrm{if}\ (c,d)=\varepsilon_{i},  \\
   0  &  \ \ \  \mathrm{otherwise}.
\end{array}\right.$$
(b) If $n\leq i$, suppose $p>n$, then for any $(c,d)\in\mathbb{Z}^{2}$ with $n\omega(c,d)\leq i$ we have
$$\ord_{p}a((\widetilde{c,d}),\varepsilon_{i})\geq i-n\omega(c,d)+n-1.$$}
\end{lem}
\textit{Proof.} The first statement follows immediately by \hyperref[lemma 2.4]{lemma 2.4(a)}, we see $(c,d)=\varepsilon_{n\omega(c,d)}$ in this setting, and $a((\widetilde{c,d}),\widetilde{\varepsilon}_{i})=\delta_{n\omega(c,d),i}$ where $\delta$ is the Kronecker delta symbol.
\par Then we show the second statement, note that $(\widetilde{c,d})\in T^{(n\omega(c,d))}$. Then by \hyperref[lemma 2.13]{lemma 2.13} we compare the coefficients of $\widetilde{\varepsilon}_{i}$ in \hyperref[2.24]{(2.24)} and obtain
$$a((\widetilde{c,d}),\widetilde{\varepsilon}_{i})=p^{i+n-n\omega(c,d)-1}\widetilde{a}((\widetilde{c,d}),\widetilde{\varepsilon}_{i})$$
with $\widetilde{a}((\widetilde{c,d}),\widetilde{\varepsilon}_{i})\in\mathbb{Z}_{q}[\widetilde{\gamma}]$, so the lemma follows. $\hfill\square$
\par With all the necessary $p$-adic estimations we need, we can study the Frobenius coefficient $\widetilde{A}(\widetilde{\varepsilon}_{j},\widetilde{\varepsilon}_{i})$ and give some decent estimations.
\begin{thm} \label{theorem 2.15}
\textit{Suppose $p>n$. \\
(a) For any $0\leq i,j\leq 2n$, $\ord_{p}\widetilde{A}(\widetilde{\varepsilon}_{j},\widetilde{\varepsilon}_{i})\geq\frac{j}{n}$.\\
(b) For any $0\leq i\leq 2n$ and $0\leq j<g(i)$, $\ord_{p}\widetilde{A}(\widetilde{\varepsilon}_{j},\widetilde{\varepsilon}_{i})>\frac{j}{n}$. In particular, if $j=g(i)$, then $\ord_{p}\widetilde{A}(\widetilde{\varepsilon}_{g(i)},\widetilde{\varepsilon}_{i})=\frac{g(i)}{n}$. \\
(c) If we furtherly restrict $p>2n^{2}-n$, then for $0\leq i\leq n$, $0\leq j\leq n-1$, we have $\ord_{p}\widetilde{A}(\widetilde{\varepsilon}_{j},\widetilde{\varepsilon}_{i})=\frac{j}{n}+\frac{(2n+1)\alpha_{i,j}}{n(p-1)}$.}
\end{thm}
\textit{Proof.} Firstly, note that $a(\widetilde{\varepsilon
}_{k},\widetilde{\varepsilon}_{k})=1$ for all $0\leq k\leq 2n$, we may rewrite \hyperref[2.25]{(2.25)} as 
\begin{equation}
\begin{aligned}
\widetilde{A}(\widetilde{\varepsilon}_{j},\widetilde{\varepsilon}_{i})=A(\widetilde{\varepsilon}_{j},\widetilde{\varepsilon}_{i})+\sum_{\begin{subarray}{l}\ \ \ \ \ (c,d)\in\mathbb{Z}^{2}\\ n\omega(c,d)\leq j, (c,d)\neq\varepsilon_{j} 
\end{subarray}}A((\widetilde{c,d}),\widetilde{\varepsilon}_{i})a((\widetilde{c,d}),\widetilde{\varepsilon}_{j})+\sum_{\begin{subarray}{l}\ \ \ (c,d)\in\mathbb{Z}^{2} \\ n\omega(c,d)\geq j+1    
\end{subarray}}A((\widetilde{c,d}),\widetilde{\varepsilon}_{i})a((\widetilde{c,d}),\widetilde{\varepsilon}_{j}).        \label{2.27}
\end{aligned}    
\end{equation}
\par Then for statement (a), by \hyperref[lemma 2.8]{lemma 2.8} we have $\ord_{p}A(\widetilde{\varepsilon}_{j},\widetilde{\varepsilon}_{i})\geq\frac{j}{n}$. Since for any $(c,d)\in\mathbb{Z}^{2}$, $a((\widetilde{c,d}),\widetilde{\varepsilon}_{j})\in\mathbb{Z}_{q}[\widetilde{\gamma}]$, apply \hyperref[lemma 2.8]{lemma 2.8} again we obtain
$$\ord_{p}A((\widetilde{c,d}),\widetilde{\varepsilon}_{i})a((\widetilde{c,d}),\widetilde{\varepsilon}_{j})\geq\ord_{p}A((\widetilde{c,d}),\widetilde{\varepsilon}_{i})\geq\omega(c,d)\geq\frac{j+1}{n}>\frac{j}{n}$$
for elements in the third summand of \hyperref[2.27]{(2.27)}. 
\par For elements in the second summand where $\omega(c,d)\leq\frac{j}{n}$, $(c,d)\neq\varepsilon_{j}$, if $0\leq j\leq n-1$, then by \hyperref[lemma 2.14]{lemma 2.14 (a)} $a((\widetilde{c,d}),\widetilde{\varepsilon}_{j})=0$, so the second summand vanishes. If $n\leq j$, by \hyperref[lemma 2.14]{lemma 2.14 (b)} and \hyperref[lemma 2.8]{lemma 2.8} we have
\begin{equation}
\begin{aligned}
\ord_{p}A((\widetilde{c,d}),\widetilde{\varepsilon}_{i})a((\widetilde{c,d}),\widetilde{\varepsilon}_{j})&\geq\omega(c,d)+j-n\omega(c,d)+n-1\geq\frac{j}{n}. \label{2.28}
\end{aligned}   
\end{equation}
Combining all the three estimations we obtain statement (a).
\par For statement (b), the third summand in \hyperref[2.27]{(2.27)} still have $p$-adic order strictly larger than $\frac{j}{n}$ like in statement (a). Since $j<g(i)$, then in the second summand, $\omega(c,d)\leq\frac{j}{n}<\frac{g(i)}{n}$, and apparently $(c,d)\neq\varepsilon_{g(i)}$. Therefore by \hyperref[theorem 2.9]{theorem 2.9} we have
$$\ord_{p}A((\widetilde{c,d}),\widetilde{\varepsilon}_{i})>\omega(c,d),\ \mathrm{and}\ \ord_{p}A(\widetilde{\varepsilon}_{j},\widetilde{\varepsilon}_{i})>\frac{j}{n}.$$
So the inequality in \hyperref[2.28]{(2.28)} is strict in case (b). We have $\ord_{p}\widetilde{A}(\widetilde{\varepsilon}_{j},\widetilde{\varepsilon}_{i})>\frac{j}{n}$. In particular, when $j=g(i)$, we have $\ord_{p}A(\widetilde{\varepsilon}_{g(i)},\widetilde{\varepsilon}_{i})=g(i)/n$ by \hyperref[theorem 2.10]{theorem 2.10}. The inequalites for the second and the third summands in \hyperref[2.27]{(2.27)} are still strict, so we have $\ord_{p}\widetilde{A}(\widetilde{\varepsilon}_{g(i)},\widetilde{\varepsilon}_{i})=g(i)/n$. This complete the proof of statement (b).
\par For statement (c), when $0\leq j\leq n-1$, the second summand vanishes, and by \hyperref[theorem 2.10]{theorem 2.10} $\ord_{p}A(\widetilde{\varepsilon}_{j},\widetilde{\varepsilon}_{i})=\frac{j}{n}+\frac{(2n+1)\alpha_{i,j}}{n(p-1)}$. Same as in the proof of statement (a), the $p$-adic order of the third summand is larger than $\frac{j+1}{n}$. Note that $0\leq\alpha_{i,j}\leq n-1$. When $p>2n^{2}-n$, we have $\frac{(2n+1)\alpha_{i,j}}{n(p-1)}<\frac{1}{n}$, therefore we get the strict $p$-adic order $\ord_{p}\widetilde{A}(\widetilde{\varepsilon}_{j},\widetilde{\varepsilon}_{i})=\frac{j}{n}+\frac{(2n+1)\alpha_{i,j}}{n(p-1)}$. $\hfill\square$

\section{Estimation of the Newton polygon and the proof of main theorems} \label{section2.3}
\par In this section, we give some results about the Newton polygons for $L$-functions of our family $\{f_{t}\}_{t\in\mathbb{F}_{q}^{*}}$. Recall
$$L(f_{t},T)^{-1}=\mathrm{det}(I-T\mathrm{Frob}_{0}^{(2)}|_{H^{2}(\Omega_{\mathcal{C}_{0}}^{\bullet})})$$
and $\mathrm{Frob}_{0}^{(2)}:H^{2}(\Omega_{\mathcal{C}_{0}}^{\bullet})\rightarrow H^{2}(\Omega_{\mathcal{C}_{0}}^{\bullet})$ can be viewed as
$$\alpha_{0}:\mathcal{C}_{0}/(D_{x}\mathcal{C}_{0}+D_{y}\mathcal{C}_{0})\rightarrow\mathcal{C}_{0}/(D_{x}\mathcal{C}_{0}+D_{y}\mathcal{C}_{0}).$$
In previous chapter, we give enough $p$-adic estimations for the matrix entries of $$\alpha_{1}:\mathcal{C}_{0}/(D_{x}\mathcal{C}_{0}+D_{y}\mathcal{C}_{0})\rightarrow\mathcal{C}_{0}/(D_{x}\mathcal{C}_{0}+D_{y}\mathcal{C}_{0})$$
w.r.t the basis $\{\widetilde{\varepsilon}_{i}\}_{0\leq i\leq 2n}$. Denote $\Gamma=(\widetilde{\varepsilon}_{0},\widetilde{\varepsilon}_{1},\cdot\cdot\cdot,\widetilde{\varepsilon}_{2n})^{\mathrm{Tr}}$ the column vector for the basis, then we rewrite \hyperref[2.23]{(2.23)} as $\alpha_{1}\Gamma=A\Gamma$ where $A=\{A_{ij}\}_{0\leq i,j\leq 2n}$ the $(2n+1)\times(2n+1)$ matrix with entries in $\mathbb{Z}_{q}[\widetilde{\gamma}]$ such that $A_{ij}=\widetilde{A}(\widetilde{\varepsilon}_{j},\widetilde{\varepsilon}_{i})$. Since $\alpha_{1}$ is $\Omega_{0}(\widetilde{\gamma})$-semilinear and $\alpha_{0}=\alpha_{1}^{a}$, we have
$$\alpha_{0}\Gamma=\alpha_{1}^{a}\Gamma=\alpha_{1}^{a-1}A\Gamma=\alpha_{1}^{a-2}A^{\sigma^{-1}}A\Gamma=\cdot\cdot\cdot=A^{(\sigma^{-1})^{a-1}}A^{(\sigma^{-1})^{a-2}}\cdot\cdot\cdot A^{\sigma^{-1}}A\Gamma$$
where $\sigma\in Gal(\Omega_{0}(\widetilde{\gamma})/\Omega_{1}(\widetilde{\gamma}))$ is the lift of Frobenius fixing $\zeta_{p}$ and $\widetilde{\gamma}$ with $\sigma^{a}=1$. Therefore the $q$-adic Newton polygon of $L(f_{t},T)^{-1}$ is the $q$-adic Newton polygon of $\det(I-A^{(\sigma^{-1})^{a-1}}A^{(\sigma^{-1})^{a-2}}\cdot\cdot\cdot A^{\sigma^{-1}}AT)$. We firstly study the $p$-adic Newton polygon of $\det(I-AT)$ which will be more straightforward to compute.
\par If we write $\det(I-AT)=1+b_{1}T+b_{2}T^{2}+\cdot\cdot\cdot b_{2n+1}T^{2n+1}$, then we will have
\begin{equation}
\begin{aligned}
b_{m}&=(-1)^{m}\sum_{\begin{subarray}{l}\ \ \ 
 0\leq u_{0}<u_{1}<\cdot\cdot\cdot<u_{m-1}\leq 2n \\ \delta\in S_{m}\  \mathrm{permuting}\ 0,1,\cdot\cdot\cdot, m-1   
\end{subarray}}\sgn(\delta)\prod_{i=0}^{m-1}A_{u_{i},u_{\delta(i)}}\\
&=(-1)^{m}\sum_{\begin{subarray}{l}\ \ \  
 0\leq u_{0}<u_{1}<\cdot\cdot\cdot<u_{m-1}\leq 2n \\ \delta\in S_{m}\  \mathrm{permuting}\ 0,1,\cdot\cdot\cdot, m-1   
\end{subarray}}\sgn(\delta)\prod_{i=0}^{m-1}\widetilde{A}(\widetilde{\varepsilon}_{u_{\delta(i)}},\widetilde{\varepsilon}_{u_{i}})  \label{3.1}
\end{aligned}
\end{equation}
where $S_{m}$ the permutation group permuting $m$ elements. We give a $p$-adic estimation for $1\leq m\leq n$. By \hyperref[theorem 2.15]{theorem 2.15 (a)}, when $p>n$,
\begin{align*}
\ord_{p}\prod_{i=0}^{m-1}\widetilde{A}(\widetilde{\varepsilon}_{u_{\delta(i)}},\widetilde{\varepsilon}_{u_{i}})&\geq \frac{u_{\delta(0)}+u_{\delta(1)}+\cdot\cdot\cdot+u_{\delta(m-1)}}{n}=\frac{u_{0}+u_{1}+\cdot\cdot\cdot+u_{m-1}}{n}\geq\frac{1}{n}\sum_{i=0}^{m-1}i=\frac{m(m-1)}{2n}.
\end{align*}
So to get a possible accurate $p$-adic estimation, we let $u_{i}=i$. Since we restrict $1\leq m\leq n$, suppose $p>2n^{2}-n$, using \hyperref[theorem 2.15]{theorem 2.15 (c)} we obtain
\begin{equation}
\begin{aligned}
\ord_{p}\prod_{i=0}^{m-1}\widetilde{A}(\widetilde{\varepsilon}_{\delta(i)},\widetilde{\varepsilon}_{i})=\sum_{i=0}^{m-1}\frac{\delta(i)}{n}+\frac{2n+1}{n(p-1)}\sum_{i=0}^{m-1}\alpha_{i,\delta(i)}=\frac{m(m-1)}{2n}+\frac{2n+1}{n(p-1)}\sum_{i=0}^{m-1}\alpha_{i,\delta(i)}
\end{aligned}
\end{equation}
for any $\delta\in S_{m}$. If we want the $p$-adic order of $b_{m}$ to be exactly the form above for some $\delta\in S_{m}$ when $1\leq m\leq n$, then we need to satisify the following two conditions:
\begin{equation}
\begin{aligned}
&(i)\ \mathrm{If}\ \{u_{0},u_{1},\cdot\cdot\cdot,u_{m-1}\}\neq\{0,1,\cdot\cdot\cdot,m-1\},\ \mathrm{then}\ \mathrm{for}\ \mathrm{any}\ \delta\in S_{m} \\ &\ord_{p}\prod_{i=0}^{m-1}\widetilde{A}(\widetilde{\varepsilon}_{u_{\delta(i)}},\widetilde{\varepsilon}_{u_{i}})>\ord_{p}\prod_{i=0}^{m-1}\widetilde{A}(\widetilde{\varepsilon}_{\delta(i)},\widetilde{\varepsilon}_{i})=\frac{m(m-1)}{2n}+\frac{2n+1}{n(p-1)}\sum_{i=0}^{m-1}\alpha_{i,\delta(i)}.    \label{3.3}
\end{aligned}    
\end{equation}
And
\begin{equation}
\begin{aligned}
&(ii)\ \ord_{p}\sum_{\delta\in S_{m}}\mathrm{sgn}(\delta)\prod_{i=0}^{m-1}\widetilde{A}(\widetilde{\varepsilon}_{\delta(i)},\widetilde{\varepsilon}_{i})=\frac{m(m-1)}{2n}+\frac{2n+1}{n(p-1)}\sum_{i=0}^{m-1}\alpha_{i,\delta'(i)} \\
&\mathrm{where}\ \delta'\in S^{0}_{m}=\{\delta\in S_{m}|\sum_{i=0}^{m-1}\alpha_{i,\delta(i)}\ \mathrm{is}\ \mathrm{minimal}\ \mathrm{among}\ \mathrm{all}\ \delta\in S_{m}\}.   \label{3.4}
\end{aligned}    
\end{equation}
\par For condition \hyperref[3.3]{(3.3)}, note that if any $u_{i}$ is replaced by $u_{i}+1$, the lower bound $p$-adic estimation of $\ord_{p}\prod_{i=0}^{m-1}\widetilde{A}(\widetilde{\varepsilon}_{u_{\delta(i)}},\widetilde{\varepsilon}_{u_{i}})$ will be $1/n$ larger by \hyperref[theorem 2.15]{theorem 2.15 (a)}. If we restrict $p>2n^{3}-n^{2}-n+1$, then 
\begin{equation*}
 \frac{2n+1}{n(p-1)}\sum_{i=0}^{m-1}\alpha_{i,\delta(i)}<\frac{1}{n}   
\end{equation*}for any $\delta\in S_{m}$ and any $1\leq m\leq n$. Therefore condition \hyperref[3.3]{(3.3)} will be satisfied by this restriction of $p$.
\par For condition \hyperref[3.4]{(3.4)}, note that it will be satisfied automatically if $\# S_{m}^{0}=1$. And this unique permutation is given by $g$ defined in \hyperref[2.20]{(2.20)} when $0\leq m\leq n$. Note that $g$ permutes non-negative integers $0,1,\cdot\cdot\cdot,m-1$ when $m=1,n,n+1,2n,2n+1$. We then have the special values of $\ord_{p}b_{m}$ :
\begin{thm}  \label{theorem 3.1}
\textit{Suppose $p>n$, then $\ord_{p}b_{m}=\frac{m(m-1)}{2n}$ when $m=1,n,n+1,2n,2n+1$.}
\end{thm}
\textit{Proof.} Note that $\alpha_{i,j}=0$ if and only if $j\equiv g(i)\ \mathrm{mod}\ n$. By \hyperref[theorem 2.15]{theorem 2.15 (b)}, for $m=1,n,n+1,2n,2n+1$, $$\ord_{p}\prod_{i=0}^{m-1}\widetilde{A}(\widetilde{\varepsilon}_{g(i)},\widetilde{\varepsilon}_{i})=\sum_{i=0}^{m-1}\frac{g(i)}{n}=\frac{m(m-1)}{2n}.$$
Here $\sum_{i=0}^{m-1}\alpha_{i,g(i)}=0$, so condition \hyperref[3.3]{(3.3)} is satisfied for those $m$'s without restricting $p>2n^{3}-n^{2}-n+1$.
And for other terms $\prod_{i=0}^{m-1}\widetilde{A}(\widetilde{\varepsilon}_{\delta(i)},\widetilde{\varepsilon}_{i})$ in \hyperref[3.4]{(3.4)} where $\delta\neq g$, there are at least one $i$ such that $\delta(i)<g(i)$, for this $i$, use \hyperref[theorem 2.15]{theorem 2.15 (b)} again we see $\ord_{p}\widetilde{A}(\widetilde{\varepsilon}_{\delta(i)},\widetilde{\varepsilon}_{i})>\frac{\delta(i)}{n}$, then
$$\ord_{p}\prod_{i=0}^{m-1}\widetilde{A}(\widetilde{\varepsilon}_{\delta(i)},\widetilde{\varepsilon}_{i})>\frac{m(m-1)}{2n}.$$
Therefore we get those strict $p$-adic orders.   $\hfill\square$
\par For general $n>1$, $2\leq m\leq n-1$, condition \hyperref[3.4]{(3.4)} becomes much more complicated. We need to study the first digits in $\widetilde{\gamma}$-adic for the sum in \hyperref[3.4]{(3.4)} and refine the estimation.
\par When $2\leq m\leq n-1$, we have $\varepsilon_{i}=x^{i}$ for $0\leq i\leq m$. Then for any $\delta\in S_{m}$, by \hyperref[2.24]{(2.24)} and \hyperref[lemma 2.14]{lemma 2.14 (a)} we have
\begin{equation}
\begin{aligned}
\widetilde{A}(\widetilde{\varepsilon}_{\delta(i)},\widetilde{\varepsilon}_{i})=\sum_{(c,d)\in\mathbb{Z}^{2}, n\omega(c,d)\geq j+1}A((\widetilde{c,d}),\widetilde{\varepsilon}_{i})a((\widetilde{c,d}),\widetilde{\varepsilon}_{\delta(i)})+A(\widetilde{\varepsilon}_{\delta(i)},\widetilde{\varepsilon}_{i}). \label{3.5}
\end{aligned}    
\end{equation}
By \hyperref[theorem 2.15]{theorem 2.15 (c)} and \hyperref[lemma 2.8]{lemma 2.8}, the right sum above has strictly larger $p$-adic order than $A(\widetilde{\varepsilon}_{\delta(i)},\widetilde{\varepsilon}_{i})$ when $p>2n^{2}-n$. And by \hyperref[proposition 2.2]{proposition 2.2 (a)}
\begin{equation}
\begin{aligned}
A(\widetilde{\varepsilon}_{\delta(i)},\widetilde{\varepsilon}_{i})&=\widetilde{\gamma}^{i-\delta(i)}B^{\sigma^{-1}}(p\delta(i)-i,0) \\
&=\widetilde{\gamma}^{i-\delta(i)}\sum_{(k,l,m)\in I(p\delta(i)-i,0)}a_{k}a_{l}a_{m}\sigma^{-k}(\Bar{t}) \\
&=\widetilde{\gamma}^{i-\delta(i)}a^{2}_{\alpha_{i,\delta(i)}}a_{\lceil\frac{p\delta(i)-i}{n}\rceil}\sigma^{-\alpha_{i,\delta(i)}}(\Bar{t})+\widetilde{\gamma}^{i-\delta(i)}\sum_{\begin{subarray}{l}
(k,l,m)\in I(p\delta(i)-i,0) \\ \ \ \ \ \  k>\alpha_{i,\delta(i)}    
\end{subarray}}a_{k}a_{l}a_{m}\sigma^{-k}(\Bar{t}) \\
&=\frac{\sigma^{-\alpha_{i,\delta(i)}}(\Bar{t})\widetilde{\gamma}^{(p-1)\delta(i)+(2n+1)\alpha_{i,\delta(i)}}}{(\alpha_{i,\delta(i)}!)^{2}(\lceil\frac{p\delta(i)-i}{n}\rceil !)}+\widetilde{\gamma}^{i-\delta(i)}\sum_{\begin{subarray}{l}
(k,l,m)\in I(p\delta(i)-i,0) \\ \ \ \ \ \  k>\alpha_{i,\delta(i)}    
\end{subarray}}a_{k}a_{l}a_{m}\sigma^{-k}(\Bar{t}).  \label{3.6}
\end{aligned}    
\end{equation}
Similar like in \hyperref[theorem 2.10]{theorem 2.10}, the term for $m=k=\alpha_{i,\delta(i)}$, $l=\lceil\frac{p\delta(i)-i}{n}\rceil$ has the smallest $p$-adic order and the right summand in \hyperref[3.6]{(3.6)} has strictly larger $p$-adic order. Combining \hyperref[3.5]{(3.5)} and \hyperref[3.6]{(3.6)} we see one way to satisfy the condition \hyperref[3.4]{(3.4)} is to show that 
\begin{equation}
\ord_{p}(\sum_{\delta\in S_{m}^{0}}\sgn(\delta)\prod_{i=0}^{m-1}\frac{1}{(\alpha_{i,\delta(i)}!)^{2}(\lceil\frac{p\delta(i)-i}{n}\rceil !)})=0 \label{3.7}
\end{equation}
under some favorable conditions.
\par To study this combinatoric sum, we follow the method in Zhu \cite{Zhu01}. We begin with the following lemma.
\begin{lem} \label{lemma 3.2}
\textit{For $0\leq i,j\leq m-1$, if we write $\alpha_{i,j}=i+1-m+\alpha_{m-1,j}+n\chi_{i,j}$ where $\chi_{i,j}=\lceil\frac{pj-i}{n}\rceil-\lceil\frac{pj-m+1}{n}\rceil$, then we have $$\chi_{i,j}=\left\{\begin{array}{cc}
0 & \mathrm{if}\ m-1-i\leq\alpha_{m-1,j} \\
1 &  \mathrm{if}\ m-1-i>\alpha_{m-1,j} 
\end{array}\right. .$$}
\end{lem}
\textit{Proof.} This follows immediately from the triangular inequality of the ceiling function.   $\hfill\square$
\par Now for any $\delta\in S_{m}$, we see
\begin{equation}
\begin{aligned}
\sum_{i=0}^{m-1}\alpha_{i,\delta(i)}&=\sum_{i=0}^{m-1}(i+1-m)+\sum_{i=0}^{m-1}\alpha_{m-1,\delta(i)}+\sum_{i=0}^{m-1}\chi_{i,\delta(i)} \\
&=\sum_{i=0}^{m-1}(i+1-m)+\sum_{i=0}^{m-1}\alpha_{m-1,i}+\sum_{i=0}^{m-1}\chi_{i,\delta(i)}.  \label{3.8}
\end{aligned}  
\end{equation}
The first and the second sum is fixed when $m$, $n$ are fixed. Only the last sum depends on the choice of $\delta\in S_{m}$. So to make $\sum_{i=0}^{m-1}\alpha_{i,\delta(i)}$ the smallest, we just need to make $\sum_{i=0}^{m-1}\chi_{i,\delta(i)}$ the smallest. The best possible choice is a $\delta\in S_{m}$ such that $\chi_{i,\delta(i)}=0$ for all $0\leq i\leq m-1$. 
\par Here we show the existence for such a $\delta$. Since for all $i\neq j$, we have $\alpha_{m-1,i}\neq\alpha_{m-1,j}$, $\{\alpha_{m-1,i}\}_{0\leq i\leq m-1}$ consists of distinct integers in $\mathbb{Z}\cap [0,n-1]$. For any $\delta\in S_{m}$, we know $\{m-1-\delta^{-1}(i)\}_{0\leq i\leq m-1}$ consists of exhausted distinct integers in $\mathbb{Z}\cap [0,m-1]$. $m<n$, then there must exist a $\delta'\in S_{m}$ such that
$$m-1-\delta'^{-1}(i)\leq\alpha_{m-1.i}\ \mathrm{for}\ \mathrm{all}\ 0\leq i\leq m-1,$$
which is just $m-1-i\leq\alpha_{m-1.\delta'(i)}\ \mathrm{for}\ \mathrm{all}\ 0\leq i\leq m-1$. By \hyperref[lemma 3.2]{lemma 3.2}, $\chi_{i,\delta'(i)}=0$ for all $0\leq i\leq m-1$. This means for this $\delta'$, we have $\delta'\in S_{m}^{0}$. Then by the fact that $\delta\in S_{m}^{0}$ if and only if $\chi_{i,\delta(i)}=0$ for all $0\leq i\leq m-1$. We have showed the following proposition:
\begin{prop} \label{proposition 3.3}
\textit{The following conditions are equivalent: \\
(a) $\delta\in S^{0}_{m}=\{\delta\in S_{m}|\sum\limits_{i=0}^{m-1}\alpha_{i,\delta(i)}\ \mathrm{is}\ \mathrm{minimal}\ \mathrm{among}\ \mathrm{all}\ \delta\in S_{m}\}$. \\
(b) $m-i-1\leq\alpha_{m-1,\delta(i)}$ for all $0\leq i\leq m-1$.\\
(c) $i+1-m+\alpha_{m-1,\delta(o)}=\alpha_{i,\delta(i)}$ for all $0\leq i\leq m-1$. \\
(d) $\lceil\frac{p\delta(i)-i}{n}\rceil=\lceil\frac{p\delta(i)-m+1}{n}\rceil$ for all $0\leq i\leq m-1$.} 
\end{prop}
\par Now we denote $$I^{0}(m)=\sum_{\delta\in S_{m}^{0}}\sgn(\delta)\prod_{i=0}^{m-1}\frac{1}{(\alpha_{i,\delta(i)}!)^{2}(\lceil\frac{p\delta(i)-i}{n}\rceil !)},\ \ U(m)=\prod_{i=0}^{m-1}(\alpha_{m-1,i})^{2}(\lceil\frac{pi-m+1}{n}\rceil !).$$
Note that $U(m)$ is a $p$-adic unit. Use \hyperlink{Proposition 3.3.}{proposition 3.3} we have
\begin{align*}
U(m)I^{0}(m)&=\sum_{\delta\in S_{m}^{0}}\sgn(\delta)\prod_{i=0}^{m-1}(\frac{\alpha_{m-1,\delta(i)} !}{\alpha_{i,\delta(i)} !})^{2} \\
&=\sum_{\delta\in S_{m}^{0}}\sgn(\delta)\prod_{i=0}^{m-1}[\frac{\alpha_{m-1,\delta(i)} !}{(i+1-m\alpha_{m-1,\delta(i)} )!}]^{2} \\
&=\sum_{\delta\in S_{m}^{0}}\sgn(\delta)\prod_{i=0}^{m-1}[\alpha_{m-1,\delta(i)}(\alpha_{m-1,\delta(i)}-1)\cdot\cdot\cdot(\alpha_{m-1,\delta(i)}-m+i+2)]^{2} \\
&=\sum_{\delta\in S_{m}}\sgn(\delta)\prod_{i=0}^{m-1}[\alpha_{m-1,\delta(i)}(\alpha_{m-1,\delta(i)}-1)\cdot\cdot\cdot(\alpha_{m-1,\delta(i)}-m+i+2)]^{2} \\
&=\sum_{\delta\in S_{m}}\sgn(\delta)\prod_{i=0}^{m-1}[\alpha_{m-1,i}(\alpha_{m-1,i}-1)\cdot\cdot\cdot(\alpha_{m-1,i}-m+\delta^{-1}(i)+2)]^{2}\\
&=\sum_{\delta\in S_{m}}\sgn(\delta)\prod_{i=0}^{m-1}[\alpha_{m-1,i}(\alpha_{m-1,i}-1)\cdot\cdot\cdot(\alpha_{m-1,i}-m+\delta(i)+2)]^{2} \\
&=\det \ V(\alpha_{m-1,0},\alpha_{m-1,1},\cdot\cdot\cdot,\alpha_{m-1,m-1})
\end{align*}
where $V$ is the Vanermonde-like matrix defined in \hyperref[assumption 1.4]{assumption 1.4}.
here we use the convention that $\alpha_{m-1,\delta(i)}(\alpha_{m-1,\delta(i)}-1)\cdot\cdot\cdot(\alpha_{m-1,\delta(i)}-m+i+2)=1$ for $i=m-1$. To show $\ord_{p}I^{0}(m)=0$, we need the above determinant to be non-zero. 
\par When $n$ and $m$ are fixed, if \hyperref[assumption 1.4]{assumption 1.4} is satisfied, then for $p>M_{n}(m)$, $p$-adic order \hyperref[3.7]{(3.7)} is zero, and then condition \hyperref[3.4]{(3.4)} will be satisfied, therefore we have proved a result of the Newton polygon: 
\begin{thm} \label{theorem 3.5}
\textit{For a fixed $n>1$, and a fixed $m$ with $2\leq m\leq n-1$, if \hyperref[assumption 1.4]{assumption 1.4} is satisfied and $p>\max\{M_{n}(m), 2n^{3}-n^{2}-n+1\}$, then we have $$\ord_{p}b_{m}=\frac{m(m-1)}{2n}+\frac{2n+1}{n(p-1)}\sum_{i=0}^{m-1}\alpha_{i,\delta(i)}$$
for any $\delta\in S^{0}_{m}=\{\delta\in S_{m}|\sum\limits_{i=0}^{m-1}\alpha_{i,\delta(i)}\ \mathrm{is}\ \mathrm{minimal}\ \mathrm{among}\ \mathrm{all}\ \delta\in S_{m}\}$.}
\end{thm}
\par Note that the bounds for the prime $p$ is not tide, and actually far from been a tide lower bound since we want the consistency of those results. One shall hope some better choice of basis for $H^{2}(\Omega_{\mathcal{C}_{0}}^{\bullet})$ will not restrict the base prime too much. 
\par In order to get the full Newton polygon, we need a functional equation of the reciprocal roots for the $L$-functions: 
\begin{lem} \label{lemma 3.6}
\textit{For any $t\in\mathbb{F}_{q}^{*}$,
$$\beta_{i}\longmapsto\frac{q^{2}}{\beta_{i}}, \ 1\leq i\leq 2n+1$$ is a one-to-one correspondence of the set of reciprocal zeros of $L(f_{t},T)^{-1}$ to the set of reciprocal zeros of $L(-f_{t},T)^{-1}$.}
\par \textit{In particular, if $n$ is odd, the correspondence comes from $L(f_{t},T)^{-1}$ to $L(f_{-t},T)^{-1}$, if $n$ is even, the correspondence comes from $L(f_{t},T)^{-1}$ to $L(-x^{n}+y+\frac{t}{xy},T)^{-1}$.}
\end{lem}
\textit{Proof:} Recall the toric exponential sums $S^{*}_{k}(f_{t})=\sum_{j=1}^{2n+1}\beta_{j}^{k}$ and we have $S^{*}(-f_{t})=\sum_{j=1}^{2n+1}\bar{\beta_{j}}^{k}$ where $\bar{\beta}_{j}$ the complex conjugacy of $\beta_{j}$, the relationship between the reciprocal roots follows \hyperref[theorem 1.2]{theorem 1.2 (ii)}.
\par Note that the substitution $x\mapsto -x$, $y\mapsto -y$ does not change the toric exponential sums, and hence does not change the $L$-functions. When $n$ is odd, after the substitution $f_{-t}$ becomes $-f_{t}$. We then see the correspondence of the set of reciprocal roots.
\par Similar situation for $n$ is even, the substitution $y\mapsto -y$ gives $$L(f_{t},T)^{-1}=L(x^{n}-y-\frac{t}{xy},T)^{-1},$$
hence the correspondence of recirocal roots passes to $L(-x^{n}+y+\frac{t}{xy},T)^{-1}$. $\hfill\square$ \\
\par We then prove one of the main result:\\
\textbf{Proof of \hyperref[theorem 1.5]{theorem 1.5}}: When $n$ is odd, note that the slopes of $0\leq i\leq n$ is independent of the choice of $t\in\mathbb{F}_{p}^{*}$ as we discussed in section 2, changing from $t$ to $-t$ does not affect the first $n$ slopes of the Newton polygon. We see result follows immediately from \hyperref[theorem 3.5]{theorem 3.5} and \hyperref[lemma 3.6]{Lemma 3.6}.
\par When $n$ is even, we need to compare the first $n$ slopes of the Newton polygon for $L(f_{t},T)^{-1}$ and that for $L(-x^{n}+y+\frac{t}{xy},T)^{-1}$. Note that for $-x^{n}+y+\frac{t}{xy}$, in \hyperref[lemma 2.3]{lemma 2.3}, $B(a,b)$ becomes
$$B(a,b)=\sum_{(k,l,m)\in I(a,b)}(-1)^{l}a_{k}a_{l}a_{m}.$$
The extra factor $(-1)^{l}$ for each terms does not affect the $p$-adic estimates. So \hyperref[lemma 2.3]{lemma 2.3} remains valid and the whole theory in section 2 remains valid for $-x^{n}+y+\frac{t}{xy}$. When we estimate the Newton polygon, the terms in the coefficient of the Frobenius action in \hyperref[3.6]{(3.6)} becomes
\begin{align*}
A(\widetilde{\varepsilon}_{\delta(i)},\widetilde{\varepsilon}_{i})&=(-1)^{\lceil\frac{p\delta(i)-i}{n}\rceil !}\cdot\frac{\sigma^{-\alpha_{i,\delta(i)}}(\Bar{t})\widetilde{\gamma}^{(p-1)\delta(i)+(2n+1)\alpha_{i,\delta(i)}}}{(\alpha_{i,\delta(i)}!)^{2}(\lceil\frac{p\delta(i)-i}{n}\rceil !)}+\widetilde{\gamma}^{i-\delta(i)}\sum_{\begin{subarray}{l}
(k,l,m)\in I(p\delta(i)-i,0) \\ \ \ \ \ \  k>\alpha_{i,\delta(i)}    
\end{subarray}}(-1)^{l}a_{k}a_{l}a_{m}\sigma^{-k}(\Bar{t}).
\end{align*}
Note that the $p$-adic estimation still remains the same. Therefore the first $n$ slopes of $L(f_{t},T)^{-1}$ and $L(-x^{n}+y+\frac{t}{xy},T)^{-1}$ are the same. By \hyperref[theorem 3.5]{theorem 3.5} and \hyperref[lemma 3.6]{Lemma 3.6} we get the result.  $\hfill\square$ 
\begin{cor} \label{corollary 3.7}
\textit{When $p>2n^{3}-n^{2}-n+1$, $f_{t}$ is ordinary if and only if $p\equiv 1\ \mathrm{mod}\ n$.}
\end{cor}
\textit{Proof:}  We firstly note that the Hodge polygon for the family $f_{t}$ is the lower convex hull of the points $\{(m,\frac{m(m-1)}{2n})\}_{0\leq m\leq 2n+1}$.
\par When $p\equiv 1\ \mathrm{mod}\ n$, by definition of $g$ in \hyperref[2.20]{(2.20)}, we see $g(i)=i$ for all $0\leq i\leq 2n$. Therefore \hyperref[theorem 2.15]{theorem 2.15 (b)} gives the triangular form of the Frobenius matrix and the result follows similarly like in \cite{AS02}.
\par When $f_{t}$ is ordinary, $i.e.$ the $q$-adic Newton polygon coincides with the Hodge polygon $\mathrm{HP}(\triangle).$ Without \hyperref[assumption 3.4]{assumption 3.4}, \hyperref[theorem 3.5]{theorem 3.5} indicates that for $0\leq m\leq n-1$,
$$\ord_{p}b_{m}\geq\frac{m(m-1)}{2n}+\frac{2n+1}{n(p-1)}\sum_{i=0}^{m-1}\alpha_{i,\delta(i)}$$
where $\delta\in S^{0}_{m}=\{\delta\in S_{m}|\sum_{i=0}^{m-1}\alpha_{i,\delta(i)}\ \mathrm{is}\ \mathrm{minimal}\ \mathrm{among}\ \mathrm{all}\ \delta\in S_{m}\}$. The ordinariness of $f_{t}$ shows that $\sum_{i=0}^{m-1}\alpha_{i,\delta(i)}=0$, therefore $\alpha_{i,\delta(i)}=0$ for all $0\leq i\leq m-1$. Suppose $p\not\equiv 1\ \mathrm{mod}\ n$, we see $g$ is not the identity permutation in $S_{n}$ that permuting $0,1,\cdot\cdot\cdot,n-1$. By the definition of $\alpha_{i,j}$ in \hyperref[1.5]{(1.5)}, when $0\leq i\leq n-1$, $\alpha_{i,j}=0$ if and only if $j=g(i)$. Then there exists an $m$, $0\leq m\leq n-1$, such that $\sum_{i=0}^{m-1}\alpha_{i,\delta(i)}>0$ when $\delta\in S_{m}^{0}$. This contradicts with $\ord_{p} b_{m}=\frac{m(m-1)}{2n}$. Therefore we must have $p\equiv 1\ \mathrm{mod}\ n$.   $\hfill\square$ 
\begin{rem}
This result can also be obtained using Wan's facial decomposition theory\cite{DW01} without the restriction of the base prime.    
\end{rem}
\par We then pass to study $q$-adic Newton polygon in more general cases for $f_{t}$ where $t\in\mathbb{F}_{q}^{*}$ with $q=p^{a}$. In general, if we have a matrix $M$ for some semilinear map over $\Omega_{0}(\widetilde{\gamma})$, the $p$-adic Newton polygon of $\det(I-MT)$ need not coincide with the $q$-adic Newton polygon of 
$\det(I-M^{(\sigma^{-1})^{a-1}}M^{(\sigma^{-1})^{a-2}}\cdot\cdot\cdot M^{\sigma^{-1}}MT)$. For an example when the two polygons do not coincide see Katz \cite{NK01} section 1.3. We impose \hyperref[assumption 1.6]{assumption 1.6} on the base prime so that we can get the $q$-adic Newton polygon: 
\begin{lem} \label{lemma 3.8}
\textit{For prime $p>n$ satisfies \hyperref[assumption 1.6]{assumption 1.6}, the splitting function in \hyperref[proposition 2.2]{proposition 2.2} satisfies
$$\ord_{p}(a_{i})=\frac{i}{p-1}\ \mathrm{for}\ 0\leq i\leq p+n-1.$$}
\end{lem}
\textit{Proof:} The case when $0\leq i\leq p-1$ follows from \hyperref[proposition 2.2]{proposition 2.2 (a)}. The coefficient of $x^{p+k}$ in the splitting function is $(\frac{k}{k!}+\frac{1}{(p+k)!})\gamma^{p+k}$ when $0\leq k\leq p-1$. When \hyperref[assumption 1.6]{assumption 1.6} is satisfied we see $\ord_{p}(\frac{k}{k!}+\frac{1}{(p+k)!})=1$, hence we get the result. $\hfill\square$ \\
\par The above lemma extends our control of the $p$-adic estimation of the Frobenius coefficients. On the chain level for coefficients $A(\widetilde{\varepsilon}_{j},\widetilde{\varepsilon}_{i})$, applying similar calculation in \hyperref[theorem 2.10]{theorem 2.10} we have the following proposition:
\begin{prop} \label{proposition 3.9}
\textit{Suppose $p>n$ and $p$ satisfies \hyperref[assumption 3.8]{assumption 3.8}, we have
$$\ord_{p}A(\widetilde{\varepsilon}_{j},\widetilde{\varepsilon}_{i})=\frac{j}{n}+\frac{2n+1}{n(p-1)}\alpha_{i,j}$$
when $0\leq i\leq n$, $0\leq j\leq 2n$, or when $n+1\leq i\leq 2n$, $n+1\leq j\leq 2n$. For $n+1\leq i\leq 2n$, $0\leq j\leq n$ we have
$$A(\widetilde{\varepsilon}_{j},\widetilde{\varepsilon}_{i})=\left\{\begin{array}{cc}
   \frac{j}{n}+\frac{2n+1}{n(p-1)}\alpha_{i,j}  &  \ \ \ \  \ \  \mathrm{if}\ j\neq g(i-n),  \\
   \frac{j}{n}+\frac{2n+1}{p-1}  &  \ \ \  \mathrm{if}\ j=g(i-n).
\end{array}\right.$$}
\end{prop}
\par On the cohomological level for coefficient $\widetilde{A}(\widetilde{\varepsilon}_{j},\widetilde{\varepsilon}_{i})$, apply \hyperref[lemma 2.14]{lemma 2.14} and \hyperref[proposition 3.10]{proposition 3.10}, then proceed with the similar calculation in \hyperref[theorem 2.15]{theorem 2.15} we have the following result: 
\begin{prop} \label{proposition 3.10}
\textit{When $p>2n^{2}-n$ and $p$ satisfies \hyperref[assumption 1.6]{assumption 1.6}, we have
$$\ord_{p}\widetilde{A}(\widetilde{\varepsilon}_{j},\widetilde{\varepsilon}_{i})=\frac{j}{n}+\frac{2n+1}{n(p-1)}\alpha_{i,j}$$
when $0\leq i\leq n$, $0\leq j\leq 2n$, or when $n+1\leq i\leq 2n$, $n+1\leq j\leq 2n$. For $n+1\leq i\leq 2n$, $0\leq j\leq n$ we have
$$\widetilde{A}(\widetilde{\varepsilon}_{j},\widetilde{\varepsilon}_{i})=\left\{\begin{array}{cc}
   \frac{j}{n}+\frac{2n+1}{n(p-1)}\alpha_{i,j}  &  \ \ \ \  \ \  \mathrm{if}\ j\neq g(i-n),  \\
   \frac{j}{n}+\frac{2n+1}{p-1}  &  \ \ \  \mathrm{if}\ j=g(i-n).
\end{array}\right.$$}
\end{prop}
\par The above proposition and \hyperref[theorem 2.15]{theorem 2.15 (c)} controls the $p$-adic order for all entries in the Frobenius matrix $A$ with entry $A_{ij}=\widetilde{A}(\widetilde{\varepsilon}_{j},\widetilde{\varepsilon}_{i})$. We then state a technical lemma that gives a condition when the $p$-adic Newton polygon and the $q$-adic Newton polygon coincides:
\begin{lem} \label{lemma 3.11}
\normalfont(\cite{Zhu02}, theorem 3.3) 
\textit{For a $m\times m$ matrix $M=(M_{ij})_{0\leq i,j\leq m-1}$ with entries in $\Omega_{0}(\widetilde{\gamma})$, denote $M^{[k]}$ the submatrix of $M$ consisting of its first $k$ rows and columns. Let
\begin{align*}
 \mu(M)&=\min\limits_{0\leq j\leq m-2}(\min\limits_{0\leq i\leq m-1}\ord_{p}M_{i,j+1}-\max\limits_{0\leq i\leq m-1}\ord_{p}M_{ij})   \\
 \eta(M)&=\max\limits_{0\leq k\leq m-2}(\ord_{p}\det M^{[k]}-\sum_{j=0}^{k-1}\min\limits_{0\leq i\leq k+1}\ord_{p}M_{ij}).
\end{align*} If $\mu(M)>m\cdot\eta(M)$, then $\mathrm{NP}_{p}\det(I-MT)=\mathrm{NP}_{q}\det(I-M^{(\sigma^{-1})^{a-1}}M^{(\sigma^{-1})^{a-2}}\cdot\cdot\cdot M^{\sigma^{-1}}MT)$.} 
\end{lem}
\par With this lemma, we can now prove the main result of the $q$-adic Newton polygon for $f_{t}$ with $t\in\mathbb{F}_{q}^{*}$:\\
\textbf{Proof of \hyperref[theorem 1.7]{theorem 1.7}}: We need to check the Frobenius matrix $A=(A_{ij})_{0\leq i,j\leq 2n}$ with $A_{ij}=\widetilde{A}(\widetilde{\varepsilon}_{j},\widetilde{\varepsilon}_{i})$ satisfies the above lemma. Applying \hyperref[proposition 3.10]{proposition 3.10} we obtain
$$\mu(A)=\min\limits_{0\leq j\leq 2n-1}(\min\limits_{0\leq i\leq 2n}\ord_{p}\widetilde{A}(\widetilde{\varepsilon}_{j+1},\widetilde{\varepsilon}_{i})-\max\limits_{0\leq i\leq 2n}\ord_{p}\widetilde{A}(\widetilde{\varepsilon}_{j},\widetilde{\varepsilon}_{i}))\geq\frac{1}{n}-\frac{2n+1}{p-1}.$$
\par By \hyperref[theorem 3.5]{theorem 3.5} and \hyperref[lemma 3.6]{lemma 3.6} we see that 
$$\det A^{[k]}=\left\{\begin{array}{cc}
   \frac{k(k-1)}{2n}+\frac{2n+1}{n(p-1)}\sum_{i=0}^{k-1}\alpha_{i,\delta(i)}  &  \ \ \ \  \ \  \mathrm{for}\ 1\leq k\leq n-1,  \\
  \frac{k(k-1)}{2n}+\frac{2n+1}{n(p-1)}\sum_{i=0}^{k-n-1}\alpha_{i,\delta(i)}  &  \ \ \  \mathrm{for}\ n+2\leq k\leq 2n-1,
\end{array}\right.$$
where $\delta\in S^{0}_{m}=\{\delta\in S_{m}|\sum_{i=0}^{m-1}\alpha_{i,\delta(i)}\ \mathrm{is}\ \mathrm{minimal}\ \mathrm{among}\ \mathrm{all}\ \delta\in S_{m}\}$. Then we obtain
$$\eta(A)=\max\limits_{0\leq k\leq 2n-2}(\ord_{p}\det A^{[k]}-\sum_{j=0}^{k-1}\min\limits_{0\leq i\leq k+1}\ord_{p}\widetilde{A}(\widetilde{\varepsilon}_{j},\widetilde{\varepsilon}_{i}))\leq \max\limits_{0\leq k\leq 2n-2}(\ord_{p}\det A^{[k]}-\frac{k(k-1)}{2n})\leq \frac{n(2n+1)}{p-1}.$$
When $p>4n^{4}+4n^{3}+3n^{2}+n+1$, $\mu(A)>(2n+1)\eta(A)$. The result follows from \hyperref[lemma 3.11]{Lemma 3.11}. $\hfill\square$

\bibliographystyle{amsplain}
\bibliography{bibfile}
\end{document}